\documentclass[12pt,reqno]{amsart}
\usepackage{amsmath}
\usepackage{amssymb}
\usepackage{verbatim}
\usepackage{mathrsfs}
\usepackage{graphicx} 
\usepackage{caption}
\usepackage{subcaption}
\usepackage{hyperref}
\usepackage{color}
\usepackage[noadjust]{cite}

\usepackage[top=0.5in,bottom=0.5in,left=0.7in,right=0.7in]{geometry}

\newtheorem{lemma}{Lemma}

\newtheorem{cor}[lemma]{Corollary}
\newtheorem{theorem}[lemma]{Theorem}
\newtheorem{remark}[lemma]{Remark}

\newtheorem{exmp}{Example}

\numberwithin{equation}{section}
\numberwithin{lemma}{section}

\newcommand{\C}{\mathbb{C}}    
\newcommand{\N}{\mathbb{N}}    
\newcommand{\R}{\mathbb{R}}    
\newcommand{\Z}{\mathbb{Z}}    

\newcommand{\wh}{\widehat}
\renewcommand{\le}{\leqslant}
\renewcommand{\ge}{\geqslant}

\newcommand{\ol}{\overline}
\newcommand{\la}{\langle}
\newcommand{\ra}{\rangle}
\newcommand{\bo}{\mathcal{O}} 

\newcommand{\eps}{\epsilon}

\newcommand{\er}{\eqref}

\newcommand{\mrE}{{\mathring{E}}}

\newcommand{\tpsi}{\tilde{\psi}}
\newcommand{\tphi}{\tilde{\phi}}

\newcommand{\tta}{\tilde{a}}
\newcommand{\ttb}{\tilde{b}}
\newcommand{\ttu}{\tilde{u}}
\newcommand{\ttm}{\tilde{m}}

\newcommand{\ttE}{\tilde{E}}

\newcommand{\bp}{ \begin{proof} }
	\newcommand{\ep}{\hfill \end{proof} }
\newcommand{\be}{ \begin{equation} }
\newcommand{\ee}{ \end{equation} }
\newcommand{\tp}{\mathsf{T}}
\newcommand{\dm}{\mathsf{M}} 
\newcommand{\dn}{\mathsf{N}} 

\newcommand{\sr}{\operatorname{sr}}  
\newcommand{\sm}{\operatorname{sm}}  
\newcommand{\vmo}{\operatorname{vm}}

\newcommand{\lpm}{\operatorname{lpm}}  

\newcommand{\td}{\pmb{\delta}}  

\newcommand{\cM}{\mathcal{M}}

\newcommand{\cG}{\mathcal{G}}

\newcommand{\cN}{\mathcal{N}}

\newcommand{\om}[1]{\omega_{#1}}

\newcommand{\pz}[1]{\pmb{0}}

\newcommand{\dR}{\mathbb{R}^d}

\newcommand{\dZ}{\mathbb{Z}^d}

\newcommand{\dLp}[1]{L_{#1}(\mathbb{R}^d)}

\newcommand{\fsupp}{\text{fsupp}}

\newcommand{\vertiii}[1]{{\left\vert\kern-0.25ex\left\vert\kern-0.25ex\left\vert #1
		\right\vert\kern-0.25ex\right\vert\kern-0.25ex\right\vert}}

\newcommand{\dlp}[1]{l_{#1}(\mathbb{Z}^d)}

\begin{document}
\title[Interpolatory Dual Framelets]
{Interpolatory Dual Framelets with a General Dilation Matrix}

\author{Ran Lu}
\address{School of Mathematics, Hohai University, Nanjing, China 211100.
\quad  {\tt rlu3@hhu.edu.cn}}

\thanks{The research of the author was supported by the National Natural Science Foundation of China under grant 12201178.
}

\makeatletter \@addtoreset{equation}{section} \makeatother

\begin{abstract}Interpolatory filters are of great interest in subdivision schemes and wavelet analysis. Due to the high-order linear-phase moment property, interpolatory refinement filters are often used to construct wavelets and framelets with high-order vanishing moments. In this paper, given a general dilation matrix $\dm$, we propose a method that allows us to construct a dual $\dm$-framelet from an arbitrary pair of $\dm$-interpolatory filters such that all framelet generators/high-pass filters (1) have the interpolatory properties; (2) have high-order vanishing moments. Our method is easy to implement, as the high-pass filters are either given in explicit formulas or can be obtained by solving specific linear systems. Motivated by constructing interpolatory dual framelets, we can further deduce a method to construct an interpolatory quasi-tight framelet from an arbitrary interpolatory filter. If, in addition, the refinement filters have symmetry, we will perform a detailed analysis of the symmetry properties that the high-pass filters can achieve. We will present several examples to demonstrate our theoretical results.

\end{abstract}


 \maketitle

\pagenumbering{arabic}


\section{Introduction}

Interpolatory wavelets have the special property that the coefficients in their expansions are the sampling values of the original signal. This special property makes interpolatory of interest in sampling theory and signal processing applications. People have been studying interpolatory wavelets extensively in the literature. To mention a few here, see \cite{d92,sb92,sb93} for construction of interpolatory wavelets from autocorrelated orthogonal refinable functions; \cite{hpz14,z14} for interpolatory polynomial wavelets; \cite{rs03,zk20} for construction of interpolatory wavelet packets; \cite{han9900,hs02} for designing interpolatory biorthogonal wavelet using coset filters.

Framelets, or wavelets frames, generalize wavelets by adding redundancy to the system. The redundancy of a framelet creates more flexibility for construction and makes it easier to achieve good properties such as high-order vanishing moments and symmetry. Compared to the well-studied interpolatory wavelets, interpolatory framelets are much less investigated. To mention some related work here: \cite{e07,k17,k22} investigate the construction of dual framelets (not necessarily interpolatory) from interpolatory refinement filters; \cite{hl20} provides a method of constructing interpolatory tight framelets from a family of filters with a prime dilation factor; the recent paper \cite{lu24} studies quincunx interpolatory quasi-tight framelets. The following question remains open: Given an arbitrary pair of refinement filters with a general dilation matrix, can we derive a dual framelet such that all generators have interpolatory properties? In this paper, we will perform a comprehensive analysis and fully answer this question.

\subsection{Backgrounds on Framelets}

Let us first review some basics of framelets. Throughout this paper, $\dm$ is a $d\times d$ dilation matrix, i.e., $\dm\in\Z^{d\times d}$, and its eigenvalues are all greater than one in modulus. For simplicity, let
$$d_{\dm}:=|\det(\dm)|.$$
Denote $ \dLp{2}$ the linear space of all square integrable functions and endow $\dLp{2}$ with the following inner product:
$$\la f,g\ra:=\int_{\dR}f(x)\ol{g(x)}dx,\qquad \forall f,g\in\dLp{2}.$$
Let $\phi,\tphi, \psi_1,\dots,\psi_s,\tpsi_1,\dots,\tpsi_s\in \dLp{2}$. We say that $\{\phi;\psi_1,\dots,\psi_s\}$ is \emph{an $\dm$-framelet} in $\dLp{2}$ if
there exist positive constants $C_1$ and $C_2$ such that
%
$$C_1\|f\|_{\dLp{2}}^2\le
\sum_{k\in \dZ} |\la f, \phi(\cdot-k)\ra|^2+\sum_{l=1}^s\sum_{j=0}^\infty \sum_{k\in \dZ}
|\la f, \psi_{l;j,k}\ra|^2\le C_2\|f\|_{\dLp{2}}^2, \quad  f\in \dLp{2},$$
where
$$\psi_{l;j,k}:=d_{\dm}^{j/2}\psi_l(\dm^j\cdot-k),\quad  l=1,\dots,s;\, j\in\N_0;\,k\in\dZ.$$
The function $\phi$ is called the \emph{scaling function} and the functions $\psi_1,\dots,\psi_s$ are called the \emph{framelet generators} of the $\dm$-framelet $\{\phi;\psi_1,\dots,\psi_s\}$. We say that $(\{\phi; \psi_1,\dots,\psi_s\},\{\tphi; \tpsi_1,\dots,\tpsi_s\})$ is \emph{a dual $\dm$-framelet} in $\dLp{2}$ if both $\{\phi; \psi_1,\dots,\psi_s\}$ and $\{\tphi; \tpsi_1,\dots,\tpsi_s\}$ are $\dm$-framelets in $\dLp{2}$ and satisfy
\be\label{qtf:expr}
f=\sum_{k\in \dZ} \la f, \phi(\cdot-k)\ra \tphi(\cdot-k)+
\sum_{l=1}^s\sum_{j=0}^\infty \sum_{k\in \dZ}
\la f, \psi_{l;j,k}\ra \tpsi_{l;j,k}, \qquad \forall f\in \dLp{2},
\ee
with the above series converging unconditionally in $\dLp{2}$.

It is well-known that the function $\phi$ in an $\dm$-framelet necessarily has the \emph{refinable structure}. Denote $\dlp{0}$ the linear space of all sequences $u=\{u(k)\}_{k\in \dZ}:\dZ\to \C$ with finitely many non-zero terms. We call every $u\in\dlp{0}$ a \emph{finitely supported filter}. A compactly supported distribution $\phi$ is called \emph{$\dm$-refinable} if there exists $a\in\dlp{0}$ such that
\be\label{ref}\phi(x)=d_{\dm}\sum_{k\in\dZ}a(k)\phi(\dm x-k),\quad\forall x\in\dR.\ee
The filter $a$ is called the \emph{$\dm$-refinement filter} of $\phi$. It is natural to characterize the above refinable relation using the Fourier transform. For $u\in\dlp{0}$, define its \emph{Fourier series} via
$$\wh{u}(\xi):=\sum_{k\in\dZ}u(k)e^{-ik\cdot \xi},\quad \xi\in\dR.$$
For any integrable function $f:\dR\to\C$ (that is, $\int_{\dR}|f(x)|dx<\infty$), its \emph{Fourier transform} is defined via
$$\wh{f}(\xi):=\int_{\dR} f(x) e^{-ix\cdot\xi} dx,\quad \xi\in \dR.$$
The definition of the Fourier transform is naturally extended to $\dLp{2}$ functions and tempered distributions. Using the Fourier transform, the refinement equation \eqref{ref} is equivalent to
\be\label{ref:f}\wh{\phi}(\dm^{\tp} \xi)=\wh{a}(\xi)\wh{\phi}(\xi),\qquad \xi\in \dR.\ee
If a filter $a\in\dlp{0}$ satisfies $\wh{a}(0)=1$, then we can define a compactly supported distribution $\phi$ via
\be\label{ref:phi}\wh{\phi}(\xi):=\prod_{j=1}^\infty\wh{a}((\dm^\tp)^{-j}\xi),\quad\forall\xi\in\dR.\ee
The function $\phi$ defined via \er{ref:phi} is called the \emph{standard $\dm$-refinable function/distribution} of the filter $a$. Except for spline refinable functions, the function $\phi$ defined as \er{ref:phi} does not have an analytic form. To construct dual framelets, we wonder whether the function $\phi$ in \er{ref:phi} belongs to $\dLp{2}$ and this can be determined by the \emph{$L_2$-smoothness exponent} of the filter $a$. To introduce the $L_2$-smoothness exponent of a filter, we need several notations:

\begin{itemize}
\item Suppose $f$ and $g$ are functions of $d$ variables which are smooth on a neighbourhood of $\xi_0\in\R^d$. By $f(\xi)=g(\xi)+\bo(\|\xi-\xi_0\|^m)$ as $\xi\to\xi_0$, we mean
$$\partial^\mu f(\xi_0)=\partial^\mu g(\xi_0),\quad\forall \mu\in\N_0^d\text{ with }|\mu|<m.$$
For $a\in \dlp{0}$, we say that $a$ has \emph{order $m$ sum rules with respect to the dilation factor $\dm$} if
\be\label{sr}\wh{a}(\xi+2\pi\omega)=\bo(\|\xi\|^m),\quad\xi\to 0,\quad\forall\omega\in\Omega_{\dm}\setminus\{0\},\ee
where $\Omega_{\dm}$ is the same as in \er{omega}. Denote $\sr(a,\dm):=m$ with $m$ being the largest non-negative integer such that \er{sr} holds.

\item For $u,v\in\dlp{0}$, define their \emph{convolution} via
$$[u*v](k):=\sum_{n\in\Z^d}u(k-n)v(n),\quad\forall k\in\Z^d.$$

\item Let $\dm$ be a $d\times d$ dilation matrix. For every $n\in\N$ and $a\in\dlp{0}$, define $a_n\in\dlp{0}$ via
\be\label{an}\wh{a_n}(\xi):=\wh{a}(\xi)\wh{a}(\dm^\tp\xi)\dots\wh{a}((\dm^\tp)^{n-1}\xi).\ee

\item For every $\mu=(\mu_1,\dots,\mu_d)^\tp\in\N_0^d$, define the filter $\nabla^\mu\td\in\dlp{0}$ via
\be\label{nab:td}\wh{\nabla^\mu\td}(\xi)=(1-e^{-i\xi_1})^{\mu_1}\dots(1-e^{-i\xi_d})^{\mu_d},\quad \forall \xi\in\R^d.\ee
\end{itemize}
Let $\dm$ be a $d\times d$ dilation matrix and suppose $\sr(a,\dm)=m$. For $1\le p\le\infty$, we define \emph{the $\ell_p$-joint spectral radius of $a$} by
\be\label{spec:a}\rho_{m}(a,\dm)_p:=|\det(\dm)|\sup\left\{ \lim_{n\to\infty}\|(\nabla^{\mu}\td)*a_n\|_{\ell_p(\Z^d)}^{1/n}:\mu\in\N_0^d, |\mu|=m\right\}.\ee
We define \emph{the $L_p$-smoothness exponent of $a$} by
\be\label{sm:a}\sm_p(a,\dm):=\begin{cases}\frac{d}{p}-\log_{\rho(\dm)}(\rho_{m}(a,\dm)_p),&1\le p<\infty,\\
-\log_{\rho(\dm)}(\rho_{m}(a,\dm)_\infty),& p=\infty.\end{cases}	\ee
where $\rho(\dm)$ is the spectral radius of $\dm$.  If $\sm_2(a,\dm)>0$, then the refinable function $\phi$ defined as in \er{ref:phi} satisfies $\|\wh{\phi}\|_{H^\tau(\R^d)}:=\left(\int_{\R^d}|\wh{\phi}(\xi)|^2(1+\|\xi\|^2)^\tau d\xi\right)^{1/2}<\infty$ whenever $0\le\tau<\sm(a,\dm)$, which in particular implies $\phi\in\dLp{2}$. See, e.g., \cite{han03,han03-1,hanbook,jj99,jj03} and many references therein for details about properties of and algorithms to compute $\sm_2(a,\dm)$.

A dual framelet is often constructed from a pair of refinement filters through \emph{extension principle}. One popular choice is the following \emph{Mixed Extension Principle} (see, e.g., \cite{dhrs03, hanbook,rs97}):

\begin{theorem}\label{thm:dft}
Let $a,\tta,b_1,\dots,b_s,\ttb_1,\dots,\ttb_s\in \dlp{0}$ be such that $\wh{a}(0)=\wh{\tta}(0)=1$. Define $\phi,\tphi,\psi_1,\dots,$ $\psi_s,\tpsi_1,\dots,\tpsi_s$ by
\be\label{ref:tphi}\wh{\phi}(\xi):=\prod_{j=1}^\infty\wh{a}((\dm^\tp)^{-j}\xi),\quad \wh{\tphi}(\xi):=\prod_{j=1}^\infty\wh{\tta}((\dm^\tp)^{-j}\xi),\ee
\be\label{ref:tpsi}\wh{\psi_l}(\xi):=\wh{b_l}(\dm^{-\tp}\xi)\wh{\phi}(\dm^{-\tp}\xi),\quad \wh{\tpsi_l}(\xi):=\wh{\ttb_l}(\dm^{-\tp}\xi)\wh{\tphi}(\dm^{-\tp}\xi),\quad l=1,\dots,s,\ee
for all $\xi\in\R^d$. Then $(\{\phi; \psi_1,\dots,\psi_s\},\{\tphi; \tpsi_1,\dots,\tpsi_s\})$ is a dual $\dm$-framelet in $\dLp{2}$ if the following conditions are satisfied:
\begin{enumerate}

\item[(1)] $\phi,\tphi\in\dLp{2}$.

\item[(2)]  $\wh{b_l}(0)=\wh{\ttb_l}(0)=0$ holds for all $l=1,\dots,s$.

\item[(3)] $(\{a;b_1,\dots,b_s\},\{\tta;\ttb_1,\dots,\ttb_s\})$ is a \emph{dual $\dm$-framelet filter bank}, i.e.,
\be \label{dffb1}
\ol{\wh{a}(\xi)}\wh{\tta}(\xi+2\pi\omega)+\sum_{l=1}^s\ol{\wh{b_l}(\xi)}\wh{\ttb_l}(\xi+2\pi\omega)=\td(\omega),\quad \xi\in\R^d,\quad\omega\in\Omega_{\dm},\ee
where
\be\label{omega}\Omega_{\dm}:=\{\omega_1,\dots,\omega_{d_{\dm}}\}:=[\dm^{-\tp}\Z^d]\cap[0,1)^d,\text{ with }\omega_1:=0,\ee
and $\td$ is the \emph{Kronecker delta function:}
\be\label{td}\td(x)=\begin{cases}1, &x=0,\\
0, &x\ne 0\end{cases}.\ee

\end{enumerate}
\end{theorem}
For a dual $\dm$-framelet filter bank $(\{a;b_1,\dots,b_s\},\{\tta;\ttb_1,\dots,\ttb_s\})$, we call $a,\tta$ the \emph{low-pass filters} and $b_1,\dots,b_s,\ttb_1,\dots,\ttb_s$ the \emph{high-pass filters}. If items (1)-(3) of Theorem~\ref{thm:dft} hold, define $\phi,\tphi$ via \er{ref:tphi} and $\psi_1,\dots,\psi_s,\tpsi_1,\dots,\tpsi_s$ via \er{ref:tpsi}, we call $(\{\phi;\psi_1,\dots,\psi_s\}, \{\tphi;\tphi_1,\dots,\tphi_s\})$ the \emph{dual $\dm$-framelet associated with the filter bank $(\{a;b_1,\dots,b_s\},$ $\{\tta;\ttb_1,\dots,\ttb_s\})$}.

\subsection{Interpolatory Filters}

Let $\dm$ be a $d\times d$ dilation matrix. A filter $a\in\dlp{0}$ is called \emph{$\dm$-interpolatory} if
\be\label{int:a}a(\dm k)=d_{\dm}^{-1}\td(k),\quad \forall k\in\Z^d.\ee
Let $a\in\dlp{0}$ be an $\dm$-interpolatory filter with $\wh{a}(0)=1$. By \cite[Corollary~5.2]{han03} or \cite[Theorem 7.3.1]{hanbook}, if $\sm_\infty(a,\dm)>0$, then the standard $\dm$-refinable function $\phi$ of $a$ must be \emph{fundamental}, that is, $\phi$ is continuous and satisfies the following interpolatory condition:
\be\label{int:phi}\phi(k)=\td(k),\quad \forall k\in\Z^d.\ee
Conversely, if the standard $\dm$-refinable fucntion $\phi$ of the filter $a$ is fundamental, then $a$ must be $\dm$-inerpolatory and $\sm_\infty(a,\dm)>0$.  Generally, there is no efficient way to compute $\sm_\infty(a,\dm)$. Nevertheless, from \cite[Theorem 3.1]{han03-1} (also see \cite[Lemma 3.1]{lu24}), we have the following estimates of $\sm_\infty(a,\dm)$:
\be\label{sm:inf:2}\sm_2(a,\dm)-\frac{d}{2}\le \sm_\infty(a,\dm)\le \sm_2(a,\dm).\ee
As discussed, the quantity $\sm_2(a,\dm)$ can be computed by several existing methods, thus yielding estimates of $\sm_\infty(a,\dm)$. Furthermore, from the second inequality in \er{sm:inf:2}, any fundamental refinable function must belong to $\dLp{2}$.

One of the essential properties of an interpolatory filter is having high-order linear-phase moments(see Lemma~\ref{lem:lpm}). This property greatly facilitates the construction of wavelets and framelets with high-order vanishing moments. Based on the high-order linear-phase moment property, several methods of constructing dual framelets with high-order vanishing moments have been established (\cite{e07,k17,k22}).

\subsection{Our Contributions and Paper Structure}

Though one can construct dual framelets with high-order vanishing moments from interpolatory filters, most existing construction methods cannot guarantee that the framelet generators (or high-pass filters) also have interpolatory properties. To be specific, let $\dm$ be a $d\times d$ dilation matrix and $a,\tta,b_1,\dots,b_s,\ttb_1,\dots,\ttb_s\in\dlp{0}$ be such that $(\{a;b_1,\dots,b_s\},$ $\{\tta;\ttb_1,\dots,\ttb_s\})$ forms a dual $\dm$-framelet filter bank, then we say that the filter bank is \emph{$\dm$-interpolatory} if
\begin{itemize}
\item $a,\tta$ are $\dm$-interpolatory and $\wh{a}(0)=\wh{\tta}(0)=1$;

\item the filters $b_l,\ttb_l$ satisfy the following interpolatory properties:
\be\label{int:b}b_l(\dm k)=\ttb_l(\dm k)=0,\quad\forall k\in\Z^d\setminus\{0\},\quad l=1,\dots,s.\ee
\end{itemize}
Given an $\dm$-interpolatory dual $\dm$-framelet filter bank $(\{a;b_1,\dots,b_s\},$ $\{\tta;\ttb_1,\dots,\ttb_s\})$. Define $\phi,\tphi,\psi_1,$ $\dots,\psi_s,\tpsi_1,\dots,\tpsi_s$ via \er{ref:tphi} and \er{ref:tpsi}. If items (1)--(3) of Theorem~\ref{thm:dft} hold, then $(\{\phi;\psi_1,\dots,\psi_s\},$ $\{\tphi;\tpsi_1,\dots,\tpsi_s\})$ is an \emph{interpolatory dual framelet} in $\dLp{2}$, that is,
\begin{itemize}
\item $\phi$ and $\tphi$ are fundamental;
\item $\psi_l$ and $\tpsi$ satisfies the following interpolatory property:
\be\label{int:psi}\psi_l(k)=\tpsi_l(k)=0,\quad\forall k\in\Z^d\setminus\{0\},\quad l=1,\dots,s.\ee
\end{itemize}
To the best of our knowledge, no general method allows us to construct an interpolatory dual framelet filter bank from an arbitrary pair of interpolatory filters. The main goal of this paper is to prove that this is always possible and the derived framelets can achieve high-order vanishing moments. The structure of this paper is organized as follows: In Section~\ref{sec:main}, we first prove the main result Theorem~\ref{thm:df:int} on constructing interpolatory dual framelets with high-order vanishing moments. Motivated by Theorem~\ref{thm:df:int}, we immediately have Corollary~\ref{cor:qtf:int} that says given any interpolatory filter, one can always construct a \emph{quasi-tight framelet} that is interpolatory and has high-order vanishing moments. In Section~\ref{sec:sym}, we perform a detailed analysis of the symmetry property of the interpolatory dual framelets. Specifically, by adding the symmetry property to the refinement filter, we will see how the high-pass filters can also achieve symmetry. Finally, in Section~\ref{sec:exmp}, we will provide several illustrative examples.

\section{Construction of Interpolatory Dual Framelets}\label{sec:main}

Given a pair $a,\tta\in\dlp{0}$ of filters such that $\wh{a}(0)=\wh{\tta}(0)=1$ and both are $\dm$-interpolatory. We study how to construct filters $b_1,\dots,b_s,\ttb_1,\dots,\ttb_s\in\dlp{0}$ such that $(\{a;b_1,\dots,b_s\},\{\tta;\ttb_1,\dots,\ttb_s\})$ is a dual $\dm$-framelet filter bank and the interpolatory condition \er{int:b} holds. Define $\Omega_{\dm}$ as in \er{omega}.  Note that $(\{a;b_1,\dots,b_s\},\{\tta;\ttb_1,\dots,\ttb_s\})$ is a dual $\dm$-framelet filter bank if and only if
\be\label{dffb2}\cM_{a,\tta}(\xi)=\sum_{l=1}^s\begin{bmatrix}\ol{\wh{b_l}(\xi+2\pi\omega_1)}\\
	\ol{\wh{b_l}(\xi+2\pi\omega_2)}\\
	\vdots\\
	\ol{\wh{b_l}(\xi+2\pi\omega_{d_{\dm}})}
\end{bmatrix}\begin{bmatrix} \wh{\ttb_l}(\xi+2\pi\omega_1) & \wh{\ttb_l}(\xi+2\pi\omega_2) & \dots & \wh{\ttb_l}(\xi+2\pi\omega_{d_{\dm}})
\end{bmatrix},\quad \xi\in\R^d,\ee
where $\cM_{a,\tta}$ is defined as
\be\label{cm:a}\cM_{a,\tta}(\xi):=I_{d_{\dm}}-\begin{bmatrix}\ol{\wh{a}(\xi+2\pi\omega_1)}\\
	\ol{\wh{a}(\xi+2\pi\omega_2)}\\
	\vdots\\
	\ol{\wh{a}(\xi+2\pi\omega_{d_{\dm}})}
\end{bmatrix}\begin{bmatrix} \wh{\tta}(\xi+2\pi\omega_1) & \wh{\tta}(\xi+2\pi\omega_2) & \dots & \wh{\tta}(\xi+2\pi\omega_{d_{\dm}})
\end{bmatrix},\quad \xi\in\R^d.\ee
For a filter $u\in\dlp{0}$ and $\gamma\in\dZ$, define the \emph{$\gamma$-coset filter of $u$ with respect to $\dm$} via
$$u^{[\gamma,\dm]}(k):=u(\gamma+\dm k),\qquad\forall k\in\dZ.$$
Using the definition of the Fourier series of $u$, it is easy to see that 
\be\label{coset:u}\wh{u}(\xi):=\sum_{\gamma\in\Gamma_{\dm}}\wh{u^{[\gamma,\dm]}}(\dm^{\tp}\xi)e^{-i\gamma\cdot\xi},\qquad\forall \xi\in\dR,\ee
where $\Gamma_{\dm}$ is a complete set of representatives of the quotient group $\dZ/[\dm\dZ]$ and is given by
\be\label{ga:dn}\Gamma_{\dm}:=[\dm[0,1)^d]\cap\dZ:=\{\gamma_1,\dots,\gamma_{d_{\dm}}\}\text{ with }\gamma_1:=0.\ee
By \er{coset:u}, we have
\be\label{coset:0}[\wh{u}(\xi+2\pi\omega_1),\, \dots,\, \wh{u}(\xi+2\pi\omega_{d_{\dm}})]=[\wh{u^{[\gamma_{1},\dm]}}(\dm^\tp\xi),\, \dots,\, \wh{u^{[\gamma_{d_{\dm}},\dm]}}(\dm^\tp\xi)]F(\xi),\quad\forall\xi\in\dR,\ee
where $F(\xi)$ is the following $d_{\dm}\times d_{\dm}$ matrix:
\be\label{Fourier}F(\xi):=[e^{-i\gamma_j\cdot(\xi+2\pi\omega_l)}]_{1\le j,l\le d_{\dm}}.\ee
Observe that $F(\xi)\ol{F(\xi)}^\tp=d_{\dm}I_{d_{\dm}}$ for all $\xi\in\dR$. Hence, \er{dffb2} is equivalent to
\be\label{dffb3}\cN_{a,\tilde{a}}(\xi)=\sum_{l=1}^s\begin{bmatrix}\ol{\wh{b_l^{[\gamma_1,\dm]}}(\xi)}\\
	\ol{\wh{b_l^{[\gamma_2,\dm]}}(\xi)}\\
	\vdots\\
	\ol{\wh{b_l^{[\gamma_{d_{\dm},\dm]}}}(\xi)}
\end{bmatrix}\begin{bmatrix} \wh{\ttb_l^{[\gamma_1,\dm]}}(\xi) & \wh{\ttb_l^{[\gamma_2,\dm]}}(\xi) & \dots & \wh{\ttb_l^{[\gamma_{d_{\dm}},\dm]}}(\xi)
\end{bmatrix},\ee
where $\cN_{a,\tta}$ is defined as
\be\label{cn:a}\cN_{a,\tta}(\xi):=d_{\dm}^{-1}I_{d_{\dm}}-\begin{bmatrix}\ol{\wh{a^{[\gamma_1,\dm]}}(\xi)}\\
	\ol{\wh{a^{[\gamma_2,\dm]}}(\xi)}\\
	\vdots\\
	\ol{\wh{a^{[\gamma_{d_{\dm},\dm]}}}(\xi)}
\end{bmatrix}\begin{bmatrix} \wh{\tta^{[\gamma_1,\dm]}}(\xi) & \wh{\tta^{[\gamma_2,\dm]}}(\xi) & \dots & \wh{\tta^{[\gamma_{d_{\dm}},\dm]}}(\xi)
\end{bmatrix},\quad \xi\in\R^d.\ee 
Using $\gamma_1=0$ and the interpolatory property of $a$ and $\tta$, we see that
$$\wh{a^{[\gamma_1,\dm]}}(\dm^\tp\xi)=\wh{\tta^{[\gamma_1,\dm]}}(\dm^\tp\xi)=d_{\dm}^{-1},\quad\forall\xi\in\dR.$$
In this case, $\cN_{a,\tta}$ becomes
\be\label{cn:a:int}\cN_{a,\tta}(\xi):=d_{\dm}^{-1}I_{d_{\dm}}-\begin{bmatrix}d_{\dm}^{-1}\\
	\ol{\wh{a^{[\gamma_2,\dm]}}(\xi)}\\
	\vdots\\
	\ol{\wh{a^{[\gamma_{d_{\dm},\dm]}}}(\xi)}
\end{bmatrix}\begin{bmatrix} d_{\dm}^{-1} & \wh{\tta^{[\gamma_2,\dm]}}(\xi) & \dots & \wh{\tta^{[\gamma_{d_{\dm}},\dm]}}(\xi)
\end{bmatrix},\quad \xi\in\R^d.\ee

\subsection{Interpolatory Dual Framelets with High-order Vanishing Moments}

Now we discuss the way to construct an interpolatory dual $\dm$-framelet filter bank $(\{a;b_1,\dots,b_s\},\{\tta;\ttb_1,\dots,\ttb_s\})$ with high-order vanishing moments from a given pair $a,\tta\in\dlp{0}$ of $\dm$-interpolatory filters with $\wh{a}(0)=\wh{\tta}(0)=1$. 

Let $(\{\phi;\psi_1,\dots,\psi_s\},\{\tphi;\tpsi_1,\dots,\tpsi_s\})$ be a dual $\dm$-framelet in $\dLp{2}$. The sparsity of the framelet expansion \er{qtf:expr} is characterized by the vanishing moments of the generators $\psi_l,\tpsi_l$, $l=1,\dots,s$.  Let $u\in\dlp{0}$, $\psi$ be a compactly supported distribution and $m\in\N_0$. We say that

\begin{itemize}
	\item $u$ has \emph{$m$ vanishing moments} if 
	\be\label{vmo}\wh{u}(\xi)=\bo(\|\xi\|^m),\quad\xi\to 0,\ee
	and denote 
	$$\vmo(u):=\sup\{m\in\N_0:\, \text{\er{vmo} holds}\};$$
	
	\item $u$ has \emph{$m$ linear-phase moments} if
	\be\label{lpm}\wh{u}(\xi)=1+\bo(\|\xi\|^m),\quad\xi\to 0,\ee
	and denote
	$$\lpm(u):=\sup\{m\in\N_0:\, \text{\er{lpm} holds}\};$$
	
	\item $\psi$ has \emph{$m$-vanishing moments} if
	\be\label{vmo:func}\wh{\psi}(\xi)=\bo(\|\xi\|^m),\quad\xi\to 0.\ee
	and denote 
	$$\vmo(\psi)=\sup\{m\in\N_0:\, \text{\er{vmo:func} holds}\}.$$
\end{itemize}

Suppose$(\{\phi;\psi_1,\dots,\psi_s\},$ $\{\tphi;\tphi_1,\dots,\tphi_s\})$ is the dual $\dm$-framelet associated with a filter bank $(\{a;b_1,\dots,b_s\},$ $\{\tta;\ttb_1,\dots,\ttb_s\})$. On one hand, using \er{dffb1} and the refinable structure of $\phi$ and $\tphi$, it is not hard to see that
\be\label{vmo:fun}\min_{1\le l\le s}\vmo(\psi_l)\le\sr(\tta,\dm),\quad \min_{1\le l\le s}\vmo(\tpsi_l)\le\sr(a,\dm).\ee
Since $\wh{a}(0)=\wh{\tta}(0)=1$, we have $\wh{\phi}(0)=\wh{\tphi}(0)=1$ and then using \er{ref:tpsi} yields 
$$\vmo(b_l)=\vmo(\psi_l),\quad \vmo(\ttb_l)=\vmo(\tpsi_l),\quad\forall l=1,\dots,s.$$
Hence, we have
\be\label{vmo:fb}\min_{1\le l\le s}\vmo(b_l)\le\sr(\tta,\dm),\quad \min_{1\le l\le s}\vmo(\ttb_l)\le\sr(a,\dm).\ee
On the other hand, define $u_{a,\tta}\in\dlp{0}$ via 
\be\label{ua}\wh{u_{a,\tta}}(\xi):=1-\ol{\wh{a}(\xi)}\wh{\tta}(\xi),\quad\xi\in\dR.\ee
It then follows from \er{dffb1} with $\omega=0$ that
\be\label{min:vmo:0}\min_{1\le l\le s}\vmo(b_l)+\min_{1\le l\le s}\vmo(\ttb_l)\le \lpm(u_{a,\tta}).\ee

Now we discuss the way to construct an interpolatory dual $\dm$-framelet filter bank $(\{a;b_1,\dots,b_s\},$ $\{\tta;\ttb_1,\dots,\ttb_s\})$ with high-order vanishing moments from a given pair $a,\tta\in\dlp{0}$ of $\dm$-interpolatory filters with $\wh{a}(0)=\wh{\tta}(0)=1$. First, we need the following key lemma on the linear-phase moment property of an interpolatory filter.

\begin{lemma}\label{lem:lpm}Let $\dm$ be a $d\times d$ dilation matrix and define $\Omega_{\dm}$ as in \er{omega}. Let $a\in\dlp{0}$ be an $\dm$-interpolatory filter with $\wh{a}(0)=1$. Then
	\be\label{int:om}\sum_{\omega\in\Omega_{\dm}}\wh{a}(\xi+2\pi\omega)=1,\quad\forall \xi\in\dR.\ee
	Moreover, if $\sr(a,\dm)=m$ for some $m\in\N_0$, then $a$ has $m$ linear-phase moments, that is,
	\be\label{a:lpm}\wh{a}(\xi)=1+\bo(\|\xi\|^m),\quad \xi\to 0.\ee
\end{lemma}

\bp On one hand, since $a$ is $\dm$-interpolatory, we must have $\wh{a^{[0,\dm]}}(\xi)=d_{\dm}^{-1}$. On the other hand, define $F(\xi)$ as in \er{Fourier} and $\Gamma_\dm=\{\gamma_1,\dots,\gamma_{\dm}\}$ as in \er{ga:dn}, using \er{coset:0} and the fact that $F(\xi)\ol{F(\xi)}=d_{\dm}I_{d_{\dm}}$ yields
\be\label{coset:2}e^{-i\gamma_j\cdot\xi}\wh{a^{[\gamma_j,\dm]}}(\dm^\tp\xi)=d_{\dm}^{-1}\sum_{l=1}^{d_{\dm}}\wh{a}(\xi+2\pi\omega_l)e^{i\gamma_l\cdot\omega_j},\quad\forall\xi\in\dR,\, j=1,\dots,d_{\dm}.\ee
In particular, as $\gamma_1=0$ and $\omega_1=0$, we have
$$1=d_{\dm}\wh{a^{[0,\dm]}}(\dm^\tp \xi)=\sum_{l=1}^{d_{\dm}}\wh{a}(\xi+2\pi\omega),\quad\forall\xi\in\dR,$$
and this proves \er{int:om}.

If $\sr(a,\dm)=m$, then $\wh{a}(\xi+2\pi\omega)=\bo(\|\xi\|^m)$ as $\xi\to 0$ for all $\omega\in\Omega_{\dm}\setminus\{0\}$. Therefore, \er{a:lpm} follows immediately from \er{int:om}.\ep

From Lemma~\ref{lem:lpm}, if $a,\tta\in\dlp{0}$ are $\dm$-interpolatory filters with $\wh{a}(0)=\wh{\tta}(0)=1$, then
\be\label{lpm:vmo}\lpm(u_{a,\tta})\ge \min\{\sr(a,\dm),\,\sr(\tta,\dm)\},\ee
where $u_{a,\tta}\in\dlp{0}$ is defined via \er{ua}. In general, we cannot tell if $\lpm(u_{a,\tta})>\min\{\sr(a,\dm),\,\sr(\tta,\dm)\}$ holds. Indeed, equality holds in \er{lpm:vmo} for many choices of $a$ and $\tta$ (for instance, see the illustrative examples in the next section). Therefore, for our main result to work for the general setting, the best results that we can achieve on the vanishing moments of the filters $b_1,\dots,b_s,\ttb_1,\dots,\ttb_s$ are \er{vmo:fb} and
$$\min_{1\le l\le s}\vmo(b_l)+\min_{1\le l\le s}\vmo(\ttb_l)\ge\min\{\sr(a,\dm),\,\sr(\tta,\dm)\}.$$

We now state and prove the main result.

\begin{theorem}\label{thm:df:int}Let $\dm$ be a $d\times d$ dilation matrix and define $\Gamma_{\dm}$ as in \er{ga:dn}. Let $a,\tta\in\dlp{0}$ be $\dm$-interpolatory filters with $\wh{a}(0)=\wh{\tta}(0)=1$. Suppose $n_1,n_2$ are positive integers such that $n_1+n_2=\min\{\sr(a,\dm),\sr(\tta,\dm)\}$. Then one can construct a dual $\dm$-framelet filter bank $(\{a;b_1,\dots,b_s\},\{\tta;\ttb_1,\dots,\ttb_s\})$ such that the interpolatory condition \er{int:b} holds and
	\be\label{vmo:b}\min_{1\le l\le s}\vmo(b_l)\ge n_1,\quad \min_{1\le l\le s}\vmo(\ttb_l)\ge n_2.\ee
	The construction steps are as follows:
	
	\begin{enumerate}
		
		\item[(S1)]For $j=2,\dots,d_{\dm}$, define
		\be\label{cos:j}\wh{h_j}(\xi):=d_{\dm}^{-1}-d_{\dm}\ol{\wh{a^{[\gamma_j,\dm]}}(\xi)}\wh{\tta^{[\gamma_j,\dm]}}(\xi),\quad \forall\xi\in\dR.\ee
		Find $u_{j,t_j},\ttu_{j,t_j}\in\dlp{0}$, $t_j=1,\dots,s_j$ for some $s_j\in\N$ such that
		\be\label{fac:h}\wh{h_j}(\xi)=\sum_{t_j=1}^{s_j}\ol{\wh{u_{j,t_j}}(\xi)}\wh{\ttu_{j,t_j}}(\xi),\quad\forall\xi\in\dR,\ee
		and
		\be\label{mom:u}\wh{u_{j,t_j}}(\xi)=\bo(\|\xi\|^{n_1}),\quad \wh{\ttu_{j,t_j}}(\xi)=\bo(\|\xi\|^{n_2}),\quad \xi\to 0,\quad t_j=1,\dots,s_j,\ee
		for all $j=2,\dots,d_{\dm}$.

		\item[(S2)]Define the filters $b_1,\ttb_1, b_j,\ttb_j, b_{j,t_j},\ttb_{j,t_j}\in\dlp{0}$ via
		\be\label{b1}\wh{b_1}(\xi):=\wh{a}(\xi)-1,\quad \wh{\ttb_1}(\xi):=1-\wh{\tta}(\xi),\ee
		\be\label{bj}\wh{b_j}(\xi):=d_{\dm}^{-\frac{1}{2}}-d_{\dm}^{\frac{1}{2}}e^{-i\gamma_j\cdot\xi}\wh{a^{[\gamma_j,\dm]}}(\dm^\tp\xi),\quad \wh{\ttb_j}(\xi):=d_{\dm}^{-\frac{1}{2}}-d_{\dm}^{\frac{1}{2}}e^{-i\gamma_j\cdot\xi}\wh{\tta^{[\gamma_j,\dm]}}(\dm^\tp\xi),\quad j=2,\dots,d_{\dm},\ee
		\be\label{bj:tj}\wh{b_{j,t_j}}(\xi):=e^{-i\gamma_j\cdot\xi}\wh{u_{j,t_j}}(\dm^\tp\xi),\quad \wh{\ttb_{j,t_j}}(\xi):=e^{-i\gamma_j\cdot\xi}\wh{\ttu_{j,t_j}}(\dm^\tp\xi),\quad t_j=1,\dots,s_j,\, j=2,\dots,d_{\dm}.\ee
		Let
		\be\label{b}\{b_1,\dots,b_s\}:=\{b_1,\dots,b_{d_{\dm}}\}\cup\{b_{j,t_j}:\, t_j=1,\dots,s_j,\, j=2,\dots,d_{\dm}\},\ee
		\be\label{ttb}\{\ttb_1,\dots,\ttb_s\}:=\{\ttb_1,\dots,\ttb_{d_{\dm}}\}\cup\{\ttb_{j,t_j}:\, t_j=1,\dots,s_j,\, j=2,\dots,d_{\dm}\}.\ee
		Then $b_1,\dots,b_s,\ttb_1,\dots,\ttb_s$ are the desired filters that satisfy all the claims.
		
	\end{enumerate}
	
	Define $\phi$ and $\tphi$ to be the standard $\dm$-refinable functions (defined as in \er{ref:tphi}) of the filters $a$ and $\tta$. If $\sm_\infty(a,\dm)>0,\sm_\infty(\tta,\dm)>0$, then by letting $\psi_l,\tpsi_l$ as in \er{ref:tpsi} for all $l=1,\dots,s$, we obtain a dual $\dm$-framelet $(\{\phi;\psi_1,\dots,\psi_s\},$ $\{\tphi;\tpsi_1,\dots,\tpsi_s\})$ such that 
	
	\begin{itemize}
		
		\item $\phi$ and $\tphi$ are both fundamental;
		
		\item $\min_{1\le l\le s}\vmo(\psi_l)\ge n_1$ and $\min_{1\le l\le s}\vmo(\tpsi_l)\ge n_2$.
	\end{itemize}

\end{theorem}

\bp  Define $\Omega_{\dm}=\{\omega_1,\dots,\om{d_{\dm}}\}$ as in \er{omega}. Using \er{coset:0}, it is easy to deduce that
\be\label{coset:1}e^{-i\gamma_j\cdot\xi}\wh{a^{[\gamma_j;\dm]}}(\dn^\tp\xi)=d_{\dm}^{-1}\sum_{l=1}^{d_{\dm}}\wh{u}(\xi+2\pi\omega_l)e^{i\gamma_l\cdot(2\pi\omega_j)},\quad\forall j=1,\dots,d_{\dm},\quad \xi\in\dR.\ee
Let $m:=\sr(a,\dm)$. By Lemma~\ref{lem:lpm}, the filter $a$ has $m$ linear-phase moments. Hence, \er{coset:1} yields
\be\label{lpm:a}e^{-i\gamma_j\cdot\xi}\wh{a^{[\gamma_j;\dm]}}(\dm^\tp\xi)=d_{\dm}^{-1}+\bo(\|\xi\|^m),\quad \xi\to 0,\quad\forall j=1,\dots,d_{\dm}.\ee
Similarly, let $\ttm:=\sr(\tta,\dm)$, we have
\be\label{lpm:tta}e^{-i\gamma_j\cdot\xi}\wh{\tta^{[\gamma_j;\dm]}}(\dm^\tp\xi)=d_{\dm}^{-1}+\bo(\|\xi\|^{\ttm}),\quad \xi\to 0,\quad\forall j=1,\dots,d_{\dm}.\ee
Therefore, we have 
$$\wh{h_j}(\xi)=\bo(\|\xi\|^{n}),\quad\xi\to 0,\quad\forall j=1,\dots,d_{\dm},$$
where $n:=\min\{m,\ttm\}$. It then follows from \cite[Lemma 5]{dhacha} that there exist $v_{j,\alpha}\in\dlp{0}$ for all $j\in\{1,\dots,d_{\dm}\}$ and $\alpha\in\N_{0,n}^d:=\{\mu=(\mu_1,\dots,\mu_d)\in\N_0^d:\,|\mu|=\mu_1+\dots+\mu_d=n\}$ such that
\be\label{hj}\wh{h_j}(\xi)=\sum_{\alpha\in\N_{0,n}^d}\wh{\nabla^\alpha\td}(\xi)\wh{v_{j,\alpha}}(\xi),\quad\forall j=1,\dots,d_{\dm}.\ee
For each $j\in\{2,\dots,d_{\dm}\}$ and $\alpha\in\N_{0,n}^d$, choose $\mu_\alpha\in\N_{0,n_1}^d,\nu_\alpha\in\N_{0,n_2}^d$ such that $\mu_\alpha+\nu_\alpha=\alpha$, and choose $v_{j,\alpha,1},v_{j,\alpha,2}\in\dlp{0}$ such that $\wh{v_{j,\alpha}}=\ol{\wh{v_{j,\alpha,1}}}\wh{v_{j,\alpha,2}}$. Define
$$\wh{u_{j,\alpha}}(\xi):=\ol{\wh{\nabla^{\mu_\alpha}\td}(\xi)}\wh{v_{j,\alpha,1}}(\xi)\quad \wh{\ttu_{j,\alpha}}(\xi):=\wh{\nabla^{\nu_\alpha}\td}(\xi)\wh{v_{j,\alpha,2}}(\xi),\quad\forall j=2,\dots,d_{\dm},\quad \alpha\in\N_{0,n}^d.$$
Clearly,
$$\wh{u_{j,\alpha}}(\xi)=\bo(\|\xi\|^{n_1}),\quad\wh{\ttu_{j,\alpha}}(\xi)=\bo(\|\xi\|^{n_2}),\quad\xi\to 0,\quad \forall j=2,\dots,d_{\dm},\quad \alpha\in\N_{0,n}^d,$$
and
$$\wh{h_j}(\xi)=\sum_{\alpha\in\N_{0,n}^d}\ol{\wh{u_{j,\alpha}}(\xi)}\wh{\ttu_{j,\alpha}}(\xi),\quad\forall j=2,\dots,d_{\dm}.$$
By letting
$$\{u_{j,1},\dots,u_{j,s_j}\}:=\{u_{j,\alpha}:\,\alpha\in\N_{0,n}^d \},\quad \{\ttu_{j,1},\dots,\ttu_{j,s_j}\}:=\{\ttu_{j,\alpha}:\,\alpha\in\N_{0,n}^d \},\quad\forall j=2,\dots,d_{\dm},$$
all requirements of (S1) are satisfied.

We then justify (S2). Observe that a filter $u\in\dlp{0}$ satisfies $u(\dm k)=0$ for all $k\in\dZ\setminus\{0\}$ if and only if $\wh{u^{[0,\dm]}}(\xi)$ is a constant function. By calculation, we have
$$\wh{b_1^{[0,\dm]}}(\xi)=\wh{a^{[0,\dm]}}(\xi)-1=1-d_{\dm}^{-1},\quad \wh{\ttb_1^{[0,\dm]}}(\xi)=1-\wh{\tta^{[0,\dm]}}(\xi)-1=1-d_{\dm}^{-1},$$
$$\wh{b_j^{[0,\dm]}}(\xi)=\wh{\ttb_j^{[0,\dm]}}(\xi)=d_{\dm}^{-\frac{1}{2}},\quad\forall j=2,\dots,d_{\dm},$$
$$\wh{b_{j,t_j}^{[0,\dm]}}(\xi)=\wh{\ttb_{j,t_j}^{[0,\dm]}}(\xi)=0,\quad\forall t_j=1,\dots,s_j,\, j=2,\dots,d_{\dm}.$$
Hence, all filters $b_1,\dots,b_s$ defined in (S2) satisfy the interpolatory property \er{int:b}. Next, using the linear-phase moments of $a$ and $\tta$, \er{lpm:a} and \er{lpm:tta}, we have $\vmo(b_j)\ge m$ and $\vmo(\ttb_j)\ge\ttm$ for all $j=1,\dots,d_{\dm}$. Furthermore, by \er{mom:u}, we have $\vmo(b_{j,t_j})\ge n_1$ and $\vmo(\ttb_{j,t_j})\ge n_2$ for all $t_j=1,\dots,s_j$, $j=2,\dots,d_{\dm}$. Therefore, \er{vmo:b} holds. We are left to prove that $(\{a;b_1,\dots,b_s\},\{\tta;\ttb_1,\dots,\ttb_s\})$ is a dual $\dm$-framelet filter bank. By calculation, the cosets of all filters in (S2) are given by the following:
$$\wh{b_1^{[\gamma_l,\dm]}}(\xi)=\wh{a^{[\gamma_l,\dm]}}(\xi),\quad \wh{\ttb_1^{[\gamma_l,\dm]}}(\xi)=-\wh{\tta^{[\gamma_l,\dm]}}(\xi),\quad\forall j=2,\dots,d_{\dm},$$
$$\wh{b_j^{[\gamma_l,\dm]}}(\xi)=-d^{\frac{1}{2}}_{\dm}\td(l-j)\wh{a^{[\gamma_l,\dm]}}(\xi),\quad \wh{\ttb_j^{[\gamma_l,\dm]}}(\xi)=-d^{\frac{1}{2}}_{\dm}\td(l-j)\wh{\tta^{[\gamma_l,\dm]}}(\xi),\quad\forall j,l=2,\dots,d_{\dm},$$
$$\wh{b_{j,t_j}^{[\gamma_l,\dm]}}(\xi)=\td(l-j)\wh{u_{j,t_j}}(\xi),\quad \wh{\ttb_{j,t_j}^{[\gamma_l,\dm]}}(\xi)=\td(l-j)\wh{\ttu_{j,t_j}}(\xi),\quad t_j=1,\dots,s_j,\, j=2,\dots,d_{\dm}.$$
Define $\cN_{a,\tta}$ as in \er{cn:a:int} and define the $d_{\dm}\times d_{\dm}$ matrix $\cN_{b,\ttb}$ via
$$\cN_{b,\ttb}(\xi):=\sum_{j=1}^{d_{\dm}}\begin{bmatrix}\ol{\wh{b_j^{[\gamma_1,\dm]}}(\xi)}\\
\vdots\\
\ol{\wh{b_j^{[\gamma_{d_{\dm}},\dm]}}(\xi)}
\end{bmatrix}\begin{bmatrix}\wh{\ttb_j^{[\gamma_1,\dm]}}(\xi) & \dots &
\wh{\ttb_j^{[\gamma_{d_{\dm}},\dm]}}(\xi)
\end{bmatrix}+\sum_{j=2}^{d_{\dm}}\sum_{t_j=1}^{s_j}\begin{bmatrix}\ol{\wh{b_{j,t_j}^{[\gamma_1,\dm]}}(\xi)}\\
\vdots\\
\ol{\wh{b_{j,t_j}^{[\gamma_{d_{\dm}},\dm]}}(\xi)}
\end{bmatrix}\begin{bmatrix}\wh{\ttb_{j,t_j}^{[\gamma_1,\dm]}}(\xi) & \dots &
\wh{\ttb_{j,t_j}^{[\gamma_{d_{\dm}},\dm]}}(\xi)
\end{bmatrix}.$$
The entries of $\cN_{b,\ttb}$ are given as follows:
\begin{align*}[\cN_{b,\ttb}(\xi)]_{1,1}&=\sum_{j=1}^{d_{\dm}}\ol{\wh{b_j^{[0,\dm]}}(\xi)}\wh{\ttb_j^{[0,\dm]}}(\xi)+\sum_{j=2}^{d_{\dm}}\sum_{t_j=1}^{s_j}\ol{\wh{b_{j,t_j}^{[0,\dm]}}(\xi)}\wh{\ttb_{j,t_j}^{[0,\dm]}}(\xi)\\
&=-(d_{\dm}^{-1}-1)^2+d_{\dm}^{-1}(d_{\dm}-1)=d_{\dm}^{-1}-d_{\dm}^{-2}=[\cN_{a,\tta}(\xi)]_{1,1},\end{align*}
\begin{align*}[\cN_{b,\ttb}(\xi)]_{l,l}&=\sum_{j=1}^{d_{\dm}}\ol{\wh{b_j^{[\gamma_l,\dm]}}(\xi)}\wh{\ttb_j^{[\gamma_l,\dm]}}(\xi)+\sum_{j=2}^{d_{\dm}}\sum_{t_j=1}^{s_j}\ol{\wh{b_{j,t_j}^{[\gamma_l,\dm]}}(\xi)}\wh{\ttb_{j,t_j}^{[\gamma_l,\dm]}}(\xi)\\
&=-\ol{\wh{a^{[\gamma_l,\dm]}}(\xi)}\wh{\tta^{[\gamma_l,\dm]}}(\xi)+d_{\dm}\ol{\wh{a^{[\gamma_l,\dm]}}(\xi)}\wh{\tta^{[\gamma_l,\dm]}}(\xi)+\sum_{t_l=1}^{s_l}\ol{\wh{u_{l,t_l}}(\xi)}\wh{\ttu_{l,t_l}}(\xi)\\
&=d_{\dm}^{-1}-\ol{\wh{a^{[\gamma_l,\dm]}}(\xi)}\wh{\tta^{[\gamma_l,\dm]}}(\xi)\\
&=[\cN_{a,\tta}(\xi)]_{l,l},\quad\forall l=2,\dots,d_{\dm},
\end{align*}
\begin{align*}[\cN_{b,\ttb}(\xi)]_{1,l}&=\sum_{j=1}^{d_{\dm}}\ol{\wh{b_j^{[0,\dm]}}(\xi)}\wh{\ttb_j^{[\gamma_l,\dm]}}(\xi)+\sum_{j=2}^{d_{\dm}}\sum_{t_j=1}^{s_j}\ol{\wh{b_{j,t_j}^{[0,\dm]}}(\xi)}\wh{\ttb_{j,t_j}^{[\gamma_l,\dm]}}(\xi)\\
&=-(d_{\dm}^{-1}-1)\wh{\tta^{[\gamma_l,\dm]}}(\xi)-\wh{\tta^{[\gamma_l,\dm]}}(\xi)\\
&=-d_{\dm}^{-1}\wh{\tta^{[\gamma_l,\dm]}}(\xi)=[\cN_{a,\tta}(\xi)]_{1,l},\quad\forall l=2,\dots,d_{\dm},\end{align*}
\begin{align*}[\cN_{b,\ttb}(\xi)]_{l,1}&=\sum_{j=1}^{d_{\dm}}\ol{\wh{b_j^{[\gamma_l,\dm]}}(\xi)}\wh{\ttb_j^{[0,\dm]}}(\xi)+\sum_{j=2}^{d_{\dm}}\sum_{t_j=1}^{s_j}\ol{\wh{b_{j,t_j}^{[\gamma_l,\dm]}}(\xi)}\wh{\ttb_{j,t_j}^{[0,\dm]}}(\xi)\\
&=(1-d_{\dm}^{-1})\ol{\wh{a^{[\gamma_l,\dm]}}(\xi)}-\ol{\wh{a^{[\gamma_l,\dm]}}(\xi)}\\
&=-d_{\dm}^{-1}\ol{\wh{a^{[\gamma_l,\dm]}}(\xi)}=[\cN_{a,\tta}(\xi)]_{l,1},\quad\forall l=2,\dots,d_{\dm},\end{align*}
\begin{align*}
[\cN_{b,\ttb}(\xi)]_{l,k}&=\sum_{j=1}^{d_{\dm}}\ol{\wh{b_j^{[\gamma_l,\dm]}}(\xi)}\wh{\ttb_j^{[\gamma_k,\dm]}}(\xi)+\sum_{j=2}^{d_{\dm}}\sum_{t_j=1}^{s_j}\ol{\wh{b_{j,t_j}^{[\gamma_l,\dm]}}(\xi)}\wh{\ttb_{j,t_j}^{[\gamma_k,\dm]}}(\xi)\\
&=\ol{\wh{b_1^{[\gamma_l,\dm]}}(\xi)}\wh{\ttb_1^{[\gamma_k,\dm]}}(\xi)=-\ol{\wh{a^{[\gamma_l,\dm]}}(\xi)}\wh{\tta^{[\gamma_k,\dm]}}(\xi)\\
&=[\cN_{a,\tta}(\xi)]_{l,k},\quad\forall l,k=2,\dots,d_{\dm},\, l\ne k.\end{align*}
Therefore, $\cN_{b,\ttb}=\cN_{a,\tta}$, which means \er{cn:a} holds with $b_1,\dots,b_s,\ttb_1,\dots,\ttb_s$ defined as in \er{b} and \er{ttb}. 

Define $\phi$ and $\tphi$ to be the standard $\dm$-refinable functions (defined as in \er{ref:tphi}) of the filters $a$ and $\tta$. If $\sm_\infty(a,\dm)>0,\sm_\infty(\tta,\dm)>0$, then $\phi$ and $\tphi$ must be fundamental and thus $\phi,\tphi\in\dLp{2}$. By letting $\psi_l,\tpsi_l$ as in \er{ref:tpsi} for all $l=1,\dots,s$, we have
$$\min_{1\le l\le s}\vmo(\psi_l)=\min_{1\le l\le s}\vmo(b_l)\ge n_1,\quad \min_{1\le l\le s}\vmo(\tpsi_l)=\min_{1\le l\le s}\vmo(\ttb_l)\ge n_2,$$ 
and it follows from Theorem~\ref{thm:dft} that $(\{\phi;\psi_1,\dots,\psi_s\},$ $\{\tphi;\tpsi_1,\dots,\tpsi_s\})$ is a dual $\dm$-framelet in $\dLp{2}$. Moreover, noting that \er{ref:tpsi} is equivalent to
$$\psi_l(x)=\sum_{k\in\dZ}b(k)\phi(\dm x-k),\quad \tpsi_l(x)=\sum_{k\in\dZ}\ttb(k)\phi(\dm x-k),\quad\forall l=1,\dots,s,\quad x\in\dR,$$
the using the interpolatory properties of $\phi$ and $\tphi$,the interpolatory property \er{int:psi} follows immediately.\ep

\begin{remark}\begin{itemize}
		
		\item From the construction procedure and the proof of Theorem~\ref{thm:df:int}, we see that the high-pass filters $b_1,\dots,b_{d_{\dm}}$ are given in explicit formulas, and the rest of the high-pass filters can be obtained by solving specific systems of (linear) equations. Hence, our method easily applies to any dimensions and dilation matrices.

		\item From the justification of (S1), we see that in the decomposition \er{hj} of $\wh{h_j}$, the number $s_j$ of terms satisfies 
		$$s_j\le \#\N_{0,n}^d=\binom{n+d-1}{d-1},\quad\forall j=2,\dots,d_{\dm},$$
		where $\#K$ denotes the number of elements in a finite set $K$ and $n:=\min\{\sr(a,\dm),\sr(\tta,\dm)\}$. Consequently, the total number $s$ of pairs of high-pass filters has the following upper bound:
		$$s\le d_{\dm}+(d_{\dm}-1)\binom{n+d-1}{d-1}.$$
		
		\item For any $u\in\dlp{0}$,
		\begin{itemize}
			\item define its \emph{filter support} by $\fsupp(u):=[s_1,t_1]\times\dots\times[s_d,t_d]$, which is the smallest $d$-dimensional interval such that $u(k)=0$ whenever $k\notin \fsupp(u)$;
			
			\item define the filter $u^\star\in\dlp{0}$ via $\wh{u^\star}(\xi):=\ol{\wh{u}(\xi)}$ for all $\xi\in\dR$. 
		\end{itemize}  
		By our construction, we have $\fsupp(b_l)\subseteq\fsupp(a)$ for all $l=1,\dots,d_{\dm}$. The supports of the rest of the filters depend on the supports of the filters $u_{j,t_j}, \ttu_{j,t_j}$ in \er{fac:h}. When we find $u_{j,t_j}, \ttu_{j,t_j}$ that satisfy \er{fac:h}, we may impose the extra constraint that
		$$\fsupp(u_{j,t_j}^\star)+\fsupp(\ttu_{j,t_j})\subseteq \fsupp(h_j)=\fsupp([a^{[\gamma_j,\dm]}]^\star)+\fsupp(\tta^{[\gamma_j,\dm]}),$$
		for all $t_j=1,\dots,s_j,$, $j=2,\dots,d_{\dm}.$ As a consequence, the supports of the filters $b_{j,t_j},\ttb_{j,t_j}$ that are defined as \er{bj:tj} satisfy
		$$\fsupp(b_{j,t_j}^\star)+\fsupp(\ttb_{j,t_j})\subseteq \fsupp(a^\star)+\fsupp(\tta),\quad\forall t_j=1,\dots,s_j,\quad j=2,\dots,d_{\dm}.$$
		This means the support sizes of all high-pass filters are reasonably controlled because they are not significantly larger than the support sizes of the low-pass filters.

	\end{itemize}
\end{remark}

\subsection{Interpolatory Quasi-tight Framelets with High-order Vanishing Moments}

Let $(\{\phi;\psi_1,\dots,\psi_s\},$ $\{\tphi;\tpsi_1,\dots,\tpsi_s\})$ be a dual $\dm$-framelet in $\dLp{2}$. If $\phi=\tphi$ and $\psi_l=\tpsi_l$ for all $l=1,\dots,s$, then $\{\phi;\psi_1,\dots,\psi_s\}$ is called a \emph{tight $\dm$-framelet} in $\dLp{2}$. Constructing tight framelets in dimension $d\ge 2$ is a well-known challenging problem as it is intrinsically related to the spectral factorization of positive-semidefinite matrices of trigonometric polynomials (see, e.g., \cite{cpss13,cpss15}). Because of this, most existing multivariate tight framelets are constructed either through tensor products of one-dimensional tight framelets (e.g., \cite{ch01,hz14}) or from specific refinable functions (e.g., \cite{hjsz18,hl20,sz22}). To circumvent the difficulties with tight framelets, in recent years, the notion of a \emph{quasi-tight framelet} has been introduced and studied in many papers such as \cite{dhacha,hl22,lu24}. Here, a \emph{quasi-tight $\dm$-framelet} is a dual $\dm$-framelet $(\{\phi;\psi_1,\dots,\psi_s\},\{\tphi;\tpsi_1,\dots,\tpsi_s\})$ such that $\phi=\tphi$ and $\psi_l=\eps\tpsi_l$ for some $\eps_l\in\{-1,1\}$ for all $l=1,\dots,s$, and is denoted by $\{\phi;\psi_1,\dots,\psi_s\}_{\eps_1,\dots,\eps_s}$. Though a quasi-tight framelet lacks certain properties, such as energy-preserving properties, it behaves almost identically to a tight framelet. However, compared to tight framelets, constructing quasi-tight framelets is much less restrictive and complicated. As demonstrated in \cite{dhacha,hl22}, a quasi-tight framelet can be derived from an arbitrary (compactly supported) refinable function, and it is much easier to achieve nice properties such as high-order vanishing moments and the balancing property on quasi-tight framelets.

Motivated by Theorem~\ref{thm:df:int}, we discuss constructing an interpolatory quasi-tight framelet from a given interpolatory refinement filter. We first recall \cite[Lemma 4.8]{lu24} (also see \cite[Lemma 6]{dhacha}), which is critical in the construction of quasi-tight framelets with high-order vanishing moments. 

\begin{lemma}\label{lem:vm}Let $A\in\dlp{0}$ be such that $\wh{A}(\xi)=\ol{\wh{A}(\xi)}$ for all $\xi\in\dR$. Suppose $\wh{A}(\xi)=\bo(\|\xi\|^{2m})$ as $\xi\to 0$ for some $m\in\N_0$, then there exist filters $u_1,\dots,u_t\in\dlp{0}$ and $\eps_1,\dots,\eps_t\in\{-1,1\}$ for some $t\in\N$ such that 
	\be\label{vm:a}\wh{A}(\xi)=\sum_{l=1}^t\eps_l|\wh{u_l}(\xi)|^2,\quad\forall\xi\in\dR,\ee
	and
	\be\label{vm:ul}\wh{u_l}(\xi)=\bo(\|\xi\|^m),\quad\xi\to0,\quad \forall l=1,\dots,t.\ee\end{lemma} 

In the original proof of \cite[Lemma 4.8]{lu24} and \cite[Lemma 6]{dhacha}, specific constructions of the filters $u_1,\dots,u_t\in\dlp{0}$ are given, which are mainly based on solving particular systems of linear equations. Moreover, an upper bound on the number $t$ is the following:
$$t\le 3\#\left(\N_{0,2m}^d\setminus [2\N_{0,m}^d]\right)-\#\left(\N_{0,2m}^d\cap [2\N_{0,m}^d]\right)=3\binom{2m+d-1}{d-1}-\binom{m+d-1}{d-1}.$$

With Lemma~\ref{lem:vm} and using the techniques in the proof of Theorem~\ref{thm:df:int}, we have the following construction of an interpolatory quasi-tight framelet with high-order vanishing moments.

\begin{cor}\label{cor:qtf:int}Let $\dm$ be a $d\times d$ dilation matrix and define $\Gamma_{\dm}$ as in \er{ga:dn}. Let $a\in\dlp{0}$ be an $\dm$-interpolatory filter with $\wh{a}(0)=1$. Suppose $\sr(a,\dm)=2m$ for some $m\in\N$, then one can construct filters $b_1,\dots,b_s\in\dlp{0}$ such that 
	\begin{enumerate}
		\item[(1)] $\{a;b_1,\dots,b_s\}_{\eps_1,\dots,\eps_s}$ forms a quasi-tight $\dm$-framelet filter bank for some $\eps_1,\dots,\eps_s\in\{-1,1\}$, that is,
		\be \label{dffb:qtf}
		\ol{\wh{a}(\xi)}\wh{a}(\xi+2\pi\omega)+\sum_{l=1}^s\eps_l\ol{\wh{b_l}(\xi)}\wh{b_l}(\xi+2\pi\omega)=\td(\omega),\quad \xi\in\R^d,\quad\omega\in\Omega_{\dm};\ee
		
		\item[(2)] $\min_{1\le l\le s}\vmo(b_l)\ge m$ and $b_l(\dm k)=0$ for all $k\in\setminus\{0\}$ and $l=1,\dots,s.$
	\end{enumerate}
	The construction steps are as follows:
	\begin{enumerate}
		
		\item[(S1)]For $j=2,\dots,d_{\dm}$, define
		$$\wh{h_j}(\xi):=d_{\dm}^{-1}-d_{\dm}|\wh{a^{[\gamma_j,\dm]}}(\xi)|^2,\quad \forall\xi\in\dR.$$
		Find $u_{j,t_j}\in\dlp{0}$ and $\eps_{j,t_j}\in\{-1,1\}$, $t_j=1,\dots,s_j$ for some $s_j\in\N$ such that
		\be\label{fac:h:qtf}\wh{h_j}(\xi)=\sum_{t_j=1}^{s_j}\eps_{j,t_j}|\wh{u_{j,t_j}}(\xi)|^2,\quad\forall\xi\in\dR,\ee
		and
		$$\wh{u_{j,t_j}}(\xi)=\bo(\|\xi\|^{m}),\quad \xi\to 0,\quad t_j=1,\dots,s_j,$$
		for all $j=2,\dots,d_{\dm}$.

		\item[(S2)]Define the filters $b_1,b_j, b_{j,t_j}\in\dlp{0}$ via
		\be\label{b1:qtf}\wh{b_1}(\xi):=1-\wh{a}(\xi),\ee
		$$\wh{b_j}(\xi):=d_{\dm}^{-\frac{1}{2}}-d_{\dm}^{\frac{1}{2}}e^{-i\gamma_j\cdot\xi}\wh{a^{[\gamma_j,\dm]}}(\dm^\tp\xi),\quad j=2,\dots,d_{\dm},$$
		$$\wh{b_{j,t_j}}(\xi):=e^{-i\gamma_j\cdot\xi}\wh{u_{j,t_j}}(\dm^\tp\xi),\quad t_j=1,\dots,s_j,\, j=2,\dots,d_{\dm},$$
		and let
		$$\eps_1=-1,\quad \eps_j=1,\quad\forall j=2,\dots,d_{\dm}.$$
		Define
		\be\label{b:qtf}\{b_1,\dots,b_s\}:=\{b_1,\dots,b_{d_{\dm}}\}\cup\{b_{j,t_j}:\, t_j=1,\dots,s_j,\, j=2,\dots,d_{\dm}\},\ee
		\be\label{eps:qtf}\{\eps_1,\dots,\eps_s\}:=\{\eps_1,\dots,\eps_{d_{\dm}}\}\cup\{\eps_{j,t_j}:\, t_j=1,\dots,s_j,\, j=2,\dots,d_{\dm}\},\ee
		then the filters $b_1,\dots,b_s$ and $\eps_1,\dots,\eps_s\in\{-1,1\}$ meet all requirements.
		
	\end{enumerate}
	Let $\phi$ be the standard $\dm$-refinable function  (defined as in \er{ref:tphi}) of the filter $a$. If $\sm_\infty(a,\dm)>0$, then by letting $\psi_l$ as in \er{ref:tpsi} for all $l=1,\dots,s$, we obtain a quasi-tight $\dm$-framelet $\{\phi;\psi_1,\dots,\psi_s\}_{\eps_1,\dots,\eps_s}$ such that 
	
	\begin{itemize}
		
		\item $\phi$ is fundamental and $\psi_1,\dots,\psi_s$ satisfy the interpolatory property
		$$\phi_l(k)=0,\quad\forall k\in\dZ\setminus\{0\},\quad l=1,\dots,s;$$
		
		\item $\min_{1\le l\le s}\vmo(\psi_l)\ge m$.
	\end{itemize}
	
\end{cor}

\section{Analysis on the Symmetry Properties of Interpolatory Dual Framelets Derived from Theorem~\ref{thm:df:int}}\label{sec:sym}

The symmetry of a framelet system is one of the most desired features in applications. Wavelet and framelet systems with symmetry produce fewer visual artifacts than non-symmetric systems. In this section, we investigate the symmetry properties of interpolatory filters and, more importantly, whether we can achieve symmetry on the interpolatory dual or quasi-tight framelets constructed from our theoretical results.

\subsection{Interpolatory Filters and Their Symmetry Properties}

The interpolatory filter $a$ can achieve \emph{symmetry} by adding extra linear constraints. The symmetry of a multidimensional filter is associated with a symmetry group. Here, by a symmetry group, we mean a finite set $\mathcal{G}$ of $d\times d$ integer matrices that form a group under matrix multiplication. Some typical examples of symmetry groups in wavelet analysis are the following:

\begin{itemize}
	\item $\mathcal{G}=\{-I_d,I_d\}$, where $I_d$ is the $d\times d$ identity matrix;
	
	\item  The \emph{full axis symmetry group}:
	\be\label{d4}D_4:=\left\{\pm I_2,\,\pm\begin{bmatrix}1 & 0\\
		0 & -1\end{bmatrix},\,\pm\begin{bmatrix}0 & 1\\
		1 & 0\end{bmatrix},\,\pm\begin{bmatrix}0 & 1\\
		-1 &0\end{bmatrix} \right\}.\ee
	This symmetry group is associated with the quadrilateral mesh in $\Z^2$.
	
	\item The \emph{hexagon symmetry group:}
	\be\label{d6}D_6:=\left\{\pm I_2,\,\pm\begin{bmatrix}0 & 1\\
		1 & 0\end{bmatrix},\,\pm\begin{bmatrix}-1 & 1\\
		0 & 1\end{bmatrix},\,\pm\begin{bmatrix}1 & 0\\
		1 & -1\end{bmatrix},\,\pm\begin{bmatrix}0 & 1\\
		-1 & 1\end{bmatrix},\, \pm\begin{bmatrix}1 & -1\\
		1 &0\end{bmatrix}\right\}.\ee
	This symmetry group is associated with the triangular mesh in $\Z^2$.
	
\end{itemize}
Let $\mathcal{G}$ be a non-trivial symmetry (i.e., $\mathcal{G}\ne\{I_d\}$) group of $d\times d$ matrices. A filter $a\in\dlp{0}$ is 
\begin{itemize}
	
	\item \emph{$\mathcal{G}$-symmetric about a point $c\in\dR$} if
	\be\label{sym:a}a(E(k-c)+c)=a(k),\qquad\forall k\in\dZ\text{ and }E\in\mathcal{G},\ee
	\item \emph{$\mathcal{G}$-anti-symmetric about a point $c\in\dR$} if
	\be\label{asym:a}a(E(k-c)+c)=-a(k),\qquad\forall k\in\dZ\text{ and }E\in\mathcal{G}.\ee
\end{itemize}
The point $c$ in \er{sym:a} and \er{asym:a} is called the \emph{symmetry centre} of the filter $a$ and must satisfy
$$(I_d-E)c\in\dZ,\quad\forall E\in\mathcal{G}.$$
If \er{sym:a} holds, then we say that $a$ has \emph{symmetry type} $(\cG,c,1)$; if \er{asym:a} holds, then $a$ has \emph{symmetry type} $(\cG,c,-1)$.

One issue we need to consider is the compatibility of the dilation matrix $\dm$ with the symmetry group $\mathcal{G}$. For a $\mathcal{G}$-symmetric filter such that $\wh{a}(0)=1$, define its standard $\dm$-refinable function $\phi$ as \er{ref:phi}. As pointed out in \cite{han02,han03-1,han04}, the symmetry of the filter $a$ does not automatically guarantee the symmetry of $\phi$. To make $\phi$ have symmetry, we need the dilation matrix $\dm$ to be \emph{compatible} with the symmetry group $\mathcal{G}$, that is, $\dm E\dm^{-1}\in\mathcal{G}$ for all $\mathcal{G}$. In this case, by \cite[Proposition 2.1]{han04},
\begin{itemize}
	\item  if $a$ is $\mathcal{G}$-symmetric about $c\in\dR$, then $\phi$ is $\mathcal{G}$-symmetric about $c_{\phi}:=(\dm-I_d)^{-1}c$, that is,
	\be\label{sym:phi}\phi(E(x-c_{\phi})+c_{\phi})=\phi(x),\qquad\forall x\in\dR\text{ and }E\in\mathcal{G}.\ee
	
	\item  if $a$ is $\mathcal{G}$-anti-symmetric about $c\in\dR$, then $\phi$ is $\mathcal{G}$-anti-symmetric about $c_{\phi}:=(\dm-I_d)^{-1}c$, that is,
	\be\label{asym:phi}\phi(E(x-c_{\phi})+c_{\phi})=-\phi(x),\qquad\forall x\in\dR\text{ and }E\in\mathcal{G}.\ee
\end{itemize}
Here are some examples of dilation matrices that are compatible with certain symmetry groups:
\begin{itemize}
	
	\item If $\dm=MI_d$ for some $M\in\N\setminus\{1\}$, then $\dm$ is compatible with any symmetry group.
	
	\item The \emph{quincunx dilation} matrices
	\be\label{quin}M_{\sqrt{2}}:=\begin{bmatrix}1 & 1\\
		1 &-1\end{bmatrix},\quad N_{\sqrt{2}}:=\begin{bmatrix}1 & -1\\
		1 & 1\end{bmatrix},\ee
	are compatible with the symmetry groups $D_4$ and $\{I_2,-I_2\}$.
	
	\item The dilation matrix 
	\be\label{sqrt3}M_{\sqrt{3}}:=\begin{bmatrix}1 & -2\\
		2 &-1\end{bmatrix}\ee
	is compatible with the symmetry groups $D_6$ and $\{I_2,-I_2\}$.
\end{itemize}

We have the following lemma regarding its symmetry type for an interpolatory filter.

\begin{lemma}\label{lem:sym} Let $\mathcal{G}$ be a non-trivial symmetry group of $d\times d$ integer matrices which satisfies the following condition: 
	\be\label{sym:cond}\text{ there exists }\mrE\in\mathcal{G}\text{ such that }I_d-\mrE\text{ is invertible.}\ee
	Let $\dm$ be a $d\times d$ dilation matrix that is compatible with $\mathcal{G}$ and $a\in\dlp{0}$ be an $\dm$-inerpolatory filter.
	Suppose $a$ has symmetry type $(\mathcal{G},c,\eps)$ for some $c\in\dR$ and $\eps\in\{1,-1\}$, then
	\begin{enumerate}
		
		\item[(1)] if $(I_d-\mrE)c\in\dm\dZ$, then $\eps=1$ and $c=0$;
		
		\item[(2)] if $(I_d-\mrE)c\notin\dm\dZ$, then $\sr(a,\dm)\le 1$. 
		
	\end{enumerate}
	
\end{lemma}

\bp If $(I_d-\mrE)c\in\dm\dZ$, using the $\dm$-interpolatory and the symmetry properties of $a$, we have
$$d_{\dm}^{-1}=a(0)=\eps a((I_d-\mrE)c)=d_{\dm}^{-1}\eps\td(\dm^{-1}(I_d-\mrE)c).$$
As $I_d-\mrE$ is invertible, the above identities hold if and only if $c=0$ and $\eps=1$. This proves item (1).

Suppose $(I_d-\mrE)c\notin\dm\dZ$. Then there exist $\gamma\in\Gamma_{\dm}\setminus\{0\}$ and $k_0\in\dZ$ such that $(I_d-\mrE)c=\gamma+\dm k_0$. Using \er{sym:a} with $E=\mrE$, we have
\be\label{sym:0:f}e^{-i[(I_d-\mrE)c]\cdot\xi}\wh{a}(\mrE^\tp\xi)=\wh{a}(\xi),\quad\forall\xi\in\dR.\ee
Then by \er{coset:u} and the $\dm$-interpolatory property of $a$, \er{sym:0:f} further becomes
\be\label{sym:cos}\begin{aligned}\wh{a}(\xi)&=d_{\dm}^{-1}e^{-i\gamma\cdot\xi}e^{-ik_0\cdot(\dm^\tp\xi)}+\sum_{\gamma'\in\Gamma_{\dm}\setminus\{0\}}e^{-i(\gamma'+\gamma)\cdot\xi}e^{-ik_0\cdot(\dm^\tp\xi)}\wh{a^{[\gamma',\dm]}}(\dm^\tp \mrE^\tp\xi)\\
	&=d_{\dm}^{-1}e^{-i\gamma\cdot\xi}e^{-ik_0\cdot(\dm^\tp\xi)}+\sum_{\gamma''\in\Gamma_{\dm}\setminus\{\gamma\}}e^{-i\gamma''\cdot\xi}e^{-i(k'+k_0)\cdot(\dm^\tp\xi)}\wh{a^{[\gamma',\dm]}}((\mrE')^\tp\dm^\tp\xi),\quad\forall\xi\in\dR,\end{aligned}\ee
where for every $\gamma'\in\Gamma_{\dm}\setminus\{0\}$, $\gamma''$ and $k'$ are the unique elements in $\Gamma_{\dm}\setminus\{\gamma\}$ and $\dZ$ such that $\gamma''+\dm k'=\gamma+\gamma'$, and $\mrE':=\dm^{-1}\mrE \dm\in\mathcal{G}$. From \er{sym:cos}, we obtain
$$\wh{a^{[\gamma,\dm]}}(\xi)=d_{\dm}^{-1}e^{-ik_0\cdot\xi},\quad\forall\xi\in\dR.$$
We prove $\sr(a,\dm)\le 1$ by contradiction. Suppose $\sr(a,\dm)\ge 2$, then by Lemma~\ref{lem:lpm}, we must have $\lpm(a)\ge 2$ and thus \er{coset:u} yields
\be\label{cos:sym}e^{-i\gamma'\cdot\xi}\wh{a^{[\gamma',\dm]}}(\dm^\tp\xi)=d_{\dm}^{-1}+\bo(\|\xi\|^2),\quad\xi\to 0,\quad\forall \gamma'\in\Gamma_{\dm}.\ee
However 
$$e^{-i\gamma\cdot\xi}\wh{a^{[\gamma,\dm]}}(\dm^\tp\xi)=d_{\dm}^{-1}e^{-i\gamma\cdot\xi}e^{-ik_0\cdot(\dm^\tp\xi)}=d_{\dm}^{-1}e^{-i(2c)\cdot\xi},\quad\forall\xi\in\dR,$$
which has at least one non-zero first-order partial derivative at $\xi=0$ and thus contradicts \er{cos:sym}. Therefore, we must have  $\sr(a,\dm)\le 1$, and this proves (2). \ep

\begin{remark}\label{rem:sym}The condition \er{sym:cond} on the symmetry group $\mathcal{G}$ is mild. Most commonly used symmetry groups in wavelet analysis and subdivision schemes, such as $\{I_d,-I_d\}$, $D_4$, $D_6$, and many of their subgroups, satisfy this condition.

	From Lemma~\ref{lem:sym}, we see that to construct an interpolatory dual $\dm$-framelet with high-order vanishing moments, filters that are symmetric about $0$ are the only meaningful choices. 
\end{remark}

\subsection{Analysis on the Symmetry Properties of the High-pass Filters}

Suppose $a,\tta$ are $\dm$-refinement filters that have $\mathcal{G}$-symmetry for some symmetry group $\mathcal{G}$. For $d\ge 2$, it is generally challenging to construct a $d$-variate dual $\dm$-framelet filter bank with all high-pass filters having symmetry. Even if this is possible, the high-pass filters generally have several different symmetry types (which is highly possible). Based on the construction algorithm of Theorem~\ref{thm:df:int}, here let us perform a detailed analysis of the symmetry properties of the high-pass filters. From the construction procedure, except for the first pair $b_1,\ttb_1$, all other high-pass filters are constructed using the coset filters of $a$ and $\tta$. Therefore, to analyze the symmetry of the high-pass filters, we must first discuss the symmetry of the coset filters $a^{[\gamma,\dm]}$ and $\tta^{[\gamma,\dm]}$ for $\gamma\in\Gamma_{\dm}$.

Following the comments in Remark~\ref{rem:sym}, here we only discuss the case when the filter $a$ is symmetric about $0$. Suppose $a\in\dlp{0}$ has symmetry type $(\mathcal{G},0,1)$ for some non-trivial symmetry group $\mathcal{G}\ne \{I_d\}$. Let $\dm$ be a $d\times d$ dilation matrix compatible with $\mathcal{G}$. Let $\gamma\in\Gamma_{\dm}$ and $c\in\dR$ be such that $(I_d-E)c\in\dZ$ for all $E\in\mathcal{G}$. Using the symmetry of $a$, we have
$$a^{[\gamma,\dm]}(E(k-c)+c)=a(\dm Ek-\dm Ec+\dm c+\gamma)=a(\tilde{E}\dm k-\tilde{E}\dm c+\dm c+\gamma)=a(\dm k-\dm c+\tilde{E}^{-1}\dm c+\tilde{E}^{-1}\gamma),$$
where $\tilde{E}:=\dm E\dm^{-1}\in\mathcal{G}$ for all $E\in\mathcal{G}$. By choosing $c:=-\dm^{-1}\gamma$, the above identities yield
$$a^{[\gamma,\dm]}(E(k-\dm^{-1}\gamma)+\dm^{-1}\gamma)=a(\dm k+\gamma)=a^{[\gamma,\dm]}(k), \quad\forall E\in\mathcal{G}\text{ with }(I_d-E)\dm^{-1}\gamma\in\dZ.$$
As a result, if there exists a non-trvial subgroup $\cG_\gamma$ of $\mathcal{G}$ such that 
\be\label{cg:ga}(I_d-E)\dm^{-1}\gamma\in\dZ,\quad\forall E\in\cG_\gamma,\ee
then $a^{[\gamma,\dm]}$ is $\cG_\gamma$-symmetric about $\dm^{-1}\gamma$. Particularly, the coset filter $a^{[0,\dm]}$ has the same symmetry type $(\mathcal{G},0,1)$ as the filter $a$ does. Here we provide some examples to illustrate the symmetry of $a^{[\gamma,\dm]}$ with $\gamma\ne 0$:
\begin{enumerate}
	\item[1.] Let $\dm=M_{\sqrt{2}}=\begin{bmatrix}1 & 1\\
	1 &-1\end{bmatrix}$ be the quincunx dilation matrix. We have $d_{\dm}=|\det(M_{\sqrt{2}})|=2$ and 
	\be\label{ga:sq2}\Gamma_{\dm}=[M_{\sqrt{2}}[0,1)^2]\cap\Z^2=\{(0,0)^\tp,e_1\}\text{ where }e_1=(1,0)^\tp.\ee
	Let $a\in(l_0(\Z^2))$ be a filter with symmetry type $(D_4,(0,0)^\tp,1)$. By calculation, $M_{\sqrt{2}}^{-1}e_1=(1/2,1/2)^\tp$, and it is straightforward to verify that $(I_2-E)(1/2,1/2)^\tp\in\Z^2$ for all $E\in D_4$. Therefore, the coset filter $a^{[e_1,M_{\sqrt{2}}]}$ has symmetry type $(D_4,(-1/2,-1/2)^\tp,1)$.
	
	\item[2.] Let $\dm=2I_2$. We have $d_{\dm}=|\det(2I_2)|=4$ and 
	\be\label{ga:2id} \Gamma_{\dm}=[2I_2[0,1)^2]\cap\Z^2=\{(0,0)^\tp,\,e_1,\,e_2,\,e_1+e_2\}\text{ where }e_1=(1,0)^\tp,\,e_2=(0,1)^\tp.\ee
	Let $a\in(l_0(\Z^2))$ be a filter with symmetry type $(D_6,(0,0)^\tp,1)$. Define the following subgroups of $D_6$:
	\be\label{sym:2id}\mathcal{G}_1:=\left\{\pm I_2,\,\pm\begin{bmatrix}1 & -1\\
		0 &-1\end{bmatrix}\right\},\quad \mathcal{G}_2:=\left\{\pm I_2,\,\pm\begin{bmatrix}1 & 0\\
		1 &-1\end{bmatrix}\right\},\quad \mathcal{G}_3:=\left\{\pm I_2,\,\pm\begin{bmatrix}0 & 1\\
		1 &0\end{bmatrix}\right\}.\ee
	By calculation, $(I_2-E)(2^{-1}e_1)\in\Z^2$ for all $E\in\mathcal{G}_1$, $(I_2-E)(2^{-1}e_2)\in\Z^2$ for all $E\in\mathcal{G}_2$, and $(I_2-E)(2^{-1}(e_1+e_2))\in\Z^2$ for all $E\in\mathcal{G}_3$. Therefore, $a^{[e_1,2I_2]}$ has symmetry type $(\mathcal{G}_1,-2^{-1}e_1,1)=(\mathcal{G}_1,-(1/2,0)^\tp,1)$, $a^{[e_2,2I_2]}$ has symmetry type $(\mathcal{G}_2,-2^{-1}e_2,1)=(\mathcal{G}_2,-(0,1/2)^\tp,1)$, and $a^{[e_1+e_2,2I_2]}$ has symmetry type $(\mathcal{G}_3,-2^{-1}(e_1+e_2),1)=(\mathcal{G}_3,-(1/2,1/2)^\tp,1)$.

	\item[3.] Let $\dm=M_{\sqrt{3}}=\begin{bmatrix} 1 &-2\\
	2 &-1\end{bmatrix}$. We have $d_{\dm}=|\det(M_{\sqrt{3}})|=3$ and 
	\be\label{ga:sq3}\Gamma_{\dm}=[M_{\sqrt{3}}[0,1)^2]\cap\Z^2=\{(0,0)^\tp,\,-e_1,\,e_2\}.\ee
	Let $a\in(l_0(\Z^2))$ be a filter with symmetry type $(D_6,(0,0)^\tp,1)$. Define the following subgroup of $D_6$:
	\be\label{sym:sq3:1}\mathcal{H}:=\left\{I_2,\,\begin{bmatrix}0 & -1\\
		1 &-1\end{bmatrix},\,\begin{bmatrix}-1 & 1\\
		-1 &0\end{bmatrix},\,\begin{bmatrix}0 & -1\\
		-1 &0\end{bmatrix},\,\begin{bmatrix}-1 & 1\\
		0 &1\end{bmatrix},\,\begin{bmatrix}1 & 0\\
		1 &-1\end{bmatrix}\right\}.\ee
	By calculation, $(I_2-E)(M_{\sqrt{3}}^{-1}e_1)\in\Z^2$ and $(I_2-E)(M_{\sqrt{3}}^{-1}e_1)\in\Z^2$ for all $E\in\mathcal{H}$. Therefore, $a^{[-e_1,M_{\sqrt{3}}]}$ has symmetry type $(\mathcal{H},M_{\sqrt{3}}^{-1}e_1,1)=(\mathcal{H},(-1/3,-2/3)^\tp,1)$, and $a^{[e_2,M_{\sqrt{3}}]}$ has symmetry type $(\mathcal{H},-M_{\sqrt{3}}^{-1}e_2,1)=(\mathcal{H},(-2/3,-1/3)^\tp,1)$.
	
\end{enumerate}

Now suppose $a,\tta\in\dlp{0}$ are $\dm$-interpolatory filters that satisfy all assumptions of Theorem~\ref{thm:df:int}. Further, assume that $\dm$ is compatible with a non-trivial symmetry group $\cG$ and both filters $a,\tta$ have symmetry type $(\cG,0,1)$. In this case, we need the following lemma to study the symmetry property of the high-pass filters constructed by Theorem~\ref{thm:df:int}.

\begin{lemma}\label{lem:sym:hp}Let $\mathcal{G}$ be a non-trivial symmetry group of $d\times d$ integer matrices and $\dm$ be a $d\times d$ dilation matrix. Suppose $u\in\dlp{0}$ has symmetry type $(\cG,c,\eps)$ for some $c\in\dR$ and $\eps\in\{-1,1\}$. Let $\gamma\in\dZ$ and define a filter $v\in\dlp{0}$ via 
	$$\wh{v}(\xi):=e^{-\gamma\cdot\xi}\wh{u}(\dm^\tp\xi),\quad\forall\xi\in\dR,$$
then	the filter $v$ has symmetry type $(\dm\cG\dm^{-1},\gamma+\dm c, \eps)$.
\end{lemma}

\bp  Note that $u$ has symmetry type $(\cG,c,\eps)$ if and only if
$$\wh{u}(E^\tp\xi)=\eps e^{i[(I_d-E)c]\cdot\xi}\wh{u}(\xi),\quad\forall\xi\in\dR,\quad E\in\cG.$$
For each $E\in\cG$, define $\ttE:=\dm E\dm^{-1}$. Then
\begin{align*}
\wh{v}(\ttE^\tp\xi)&=e^{-i\gamma\cdot(\ttE^\tp\xi)}\wh{u}(\dm^\tp\ttE^\tp \xi)=e^{-i(\ttE\gamma)\cdot\xi}\wh{u}(E^\tp\dm^\tp \xi)=\eps e^{-i(\ttE\gamma)\cdot\xi}e^{i[(I_d-E)c]\cdot(\dm^\tp\xi)}\wh{u}(\dm^\tp \xi)\\
&=\eps e^{i[(I_d-\ttE)\gamma]\cdot\xi}e^{i[(\dm-\dm E)c]\cdot\xi}e^{-i\gamma\cdot\xi}\wh{u}(\dm^\tp\xi)=\eps e^{i[(I_d-\ttE)\gamma]\cdot\xi}e^{i[(\dm-\ttE\dm) c]\cdot\xi}\wh{v}(\xi)\\
&=\eps e^{i[(I_d-\ttE)(\gamma+\dm c)]\cdot\xi}\wh{v}(\xi),
\end{align*}
which means $v$ has symmetry type $(\dm\cG\dm^{-1},\gamma+\dm c, \eps)$.\ep

For each $\gamma\in\Gamma_{\dm}\setminus\{0\}$, assume there exists a non-trivial subgroup $\cG_{\gamma}$ of $\cG$ such that \er{cg:ga} holds. Then the coset filters $a^{[\gamma,\dm]}$ and $\tta^{[\gamma,\dm]}$ both have symmetry type $(\cG_\gamma,-\dm^{-1}\gamma,1)$. Define $b_j,\ttb_j\in\dlp{0}$, $j=2,\dots,d_{\dm}$ as in \er{bj}. If $\dm$ is compatible with $\cG$, then it follows from Lemma~\ref{lem:sym:hp} that
\be\label{sym:bj}b_j,\,\ttb_j\text{ have symmetry type }(\dm\cG_{\gamma_j}\dm^{-1},0,1),\quad\forall j=2,\dots,d_{\dm}.\ee 
Next, define $h_j\in\dlp{0}$, $j=2,\dots,d_{\dm}$ as in \er{hj}. It is trivial that $h_j$ has symmetry type $(\cG_{\gamma_j},0,1)$ for all $j=2,\dots,d_{\dm}$. As $h_j$ has symmetry, we can impose symmetry constraints to the filters $u_{j,t_j}$ and $\ttu_{j,t_j}$ that appear in the decomposition \er{fac:h}. Once we can find symmetric (or anti-symmetric) filters $u_{j,t_j},\ttu_{j,t_j}$ that meet all requirements in (S2) of Theorem~\ref{thm:df:int}, then the high-pass filters $b_{j,t_j},\ttb_{j,t_j}$ defined as in \er{bj:tj} are likely to have symmetry as well. Indeed, suppose $u_{j,t_j}$ has symmetry type $(\cG_{j,t_j},c_{j,t_j},\eps_{j,t_j},)$ and $\ttu_{j,t_j}$ has symmetry type $(\tilde{\cG}_{j,t_j},\tilde{c}_{j,t_j},\tilde{\eps}_{j,t_j})$ for some non-trivial subgroups $\cG_{j,t_j},\tilde{\cG}_{j,t_j}$ of $\cG$, $c_{j,t_j},\tilde{c}_{j,t_j}\in\dR$, and $\eps_{j,t_j},\tilde{\eps}_{j,t_j}\in\{-1,1\}$. Then if $\dm$ is compatible with $\cG$, we conclude from Lemma~\ref{lem:sym:hp} that
\be\label{sym:bj:tj}b_{j,t_j}\text{ has symmetry type }(\dm\cG_{j,t_j}\dm^{-1},\gamma_j+\dm c_{j,t_j},\eps_{j,t_j}),\quad\forall  t_j=1,\dots,s_j,\quad j=2,\dots,d_{\dm};\ee 
\be\label{sym:tbj:tj}\ttb_{j,t_j}\text{ has symmetry type }(\dm\tilde{\cG}_{j,t_j}\dm^{-1},\gamma_j+\dm \tilde{c}_{j,t_j},\tilde{\eps}_{j,t_j}),\quad\forall t_j=1,\dots,s_j,\quad j=2,\dots,d_{\dm}.\ee 

In summary, given a pair of $\dm$-interpolatory filters $a,\tta\in\dlp{0}$ that meet all assumptions of Theorem~\ref{thm:df:int} and are both $\cG$-symmetric about $(0,0)^\tp$ for some symmetry group $\cG$. Then, the high-pass filters of an interpolatory dual framelet derived through Theorem~\ref{thm:df:int} can achieve symmetry. Similarly, given an $\dm$-interpolatory filter $a\in\dlp{0}$ that is symmetric about $(0,0)^\tp$ and has $2m$ sum rules with respect to $\dm$, an interpolatory quasi-tight framelet that is constructed through Corollary~\ref{cor:qtf:int} can also achieve symmetry.

\section{Illustrative Examples for $d=2$}\label{sec:exmp}

This section presents some two-dimensional interpolatory dual framelets to illustrate Theorem~\ref{thm:df:int} and Corollary~\ref{cor:qtf:int}. To construct an interpolatory dual or quasi-tight framelet filter bank, we first need an interpolatory filter $a\in\dlp{0}$ with $\wh{a}(0)=1$. For a general $d\times d$ dilation matrix $\dm$, we can take the following steps to construct an $\dm$-interpolatory refinement filter $a$ that has order $m$ sum rules:

\begin{enumerate}
	
	\item[(S1)]Parametrize a filter $a$ by
	$$\wh{a}(\xi)=\sum_{k\in K}a(k)e^{-ik\cdot\xi},\qquad\forall \xi\in\dR,$$
	where $K=[s_1,t_1]\times\dots\times[s_d,t_d]$ for some $s_1,t_1,\dots,s_d,t_d\in\Z$. Solve the linear equation $\sum_{k\in  K}a(k)=1$ and update the filter $a$ by substituting in the solutions of the above system. If symmetry is required, then add the additional linear constraint \er{sym:a} or \er{asym:a} to obtain a filter $a$ with symmetry type $(\mathcal{G},c,\eps)$.
	
	\item[(S2)] \textbf{Sum rule conditions for $a$:} Let $m\in\N$, solve the following linear system:
	$$\partial^\mu\wh{a}(2\pi\omega)=0,\quad  \forall\omega\in\Omega_{\dm}\setminus\{0\},\,\mu\in\N_0^d\text{ with }|\mu|<m,$$
	where $\Omega_{\dm}$ is defined as in \er{omega}. Update the filter $a$ by substituting the solutions of the above system.
	
	\item[(S3)] \textbf{Interpolatory condition:} Solve the following system of linear equations:
	$$a(\dm k)=d_{\dm}^{-1}\td(k),\quad \forall k\in [\dm^{-1}K]\cap\dZ.$$
	Update the filter $a$ by substituting the solutions of the above system.
	
	\item[(S4)] \textbf{Try to optimize the smoothness exponent:} Select parameter values among the remaining free parameters such that $\sm_2(a,\dm)$ is as large as possible. Ideally, try to achieve $\sm_2(a,\dm)>\frac{d}{2}$ so that $\sm_\infty(a,\dm)>0$. If this is not possible, try to directly estimate $\sm_\infty(a,\dm)$ by using the structural properties of the filter $a$.
	
\end{enumerate}

Here, we present some examples of two-dimensional interpolatory dual framelets. For $u\in l_0(\Z^2)$, denote its \emph{filter support} by $\fsupp(u):=[s_1,s_2]\times[t_1,t_2]$ for some $s_1,s_2,t_1,t_2\in\Z$, which is the smallest two-dimensional interval such that $u(k)=0$ whenever $k\notin\fsupp(u)$. We use the following way to present a finitely supported filter $u\in l_0(\Z^2)$: suppose $\fsupp(u)=[s_1,s_2]\times[t_1,t_2]$, then we write
$$u=\begin{bmatrix}u(s_1,t_2) & u(s_1+1,t_2) &  \dots & u(k_2,t_2)\\
u(s_1,t_2-1) & u(s_1+1,t_2-1) &  \dots & u(s_2,t_2-1)\\
\vdots & \vdots &\ddots &\vdots\\
u(s_1,t_1) & u(s_1+1,t_1) &  \dots & u(s_2,t_1)\end{bmatrix}_{[s_1,s_2]\times[t_1,t_2]}.$$
For example, $\wh{u}(\xi_1,\xi_2)=e^{-i\xi_1}-e^{i\xi_2}$ is presented as $u=\begin{bmatrix}
0 & 1\\
-1 &0
\end{bmatrix}_{[0,1]\times[-1,0]}$.

\begin{exmp}\label{ex1}Let $\dm=M_{\sqrt{2}}=\begin{bmatrix}1 & 1\\
	1 &-1\end{bmatrix}$ be the quincunx dilation matrix. Parametrize an $M_{\sqrt{2}}$-ineroplatory filter $A\in l_0(\Z^2)$ such that
	\begin{itemize}
		\item $\wh{A}(0,0)=1$ and $\sr(A,M_{\sqrt{2}})=4$; 
		
		\item $A$ is $D_4$-symmetric about $(0,0)^\tp$;
	\end{itemize}
	As the following:
	\be\label{A:ex1}{\footnotesize A=\begin{bmatrix}0 & t_1 & 0 & t_2 & 0 & t_1 & 0\\[0.2cm]
			t_1 & 0 &  -\frac{1}{64}-3t_1-t_2 &  0 &  -\frac{1}{64}-3t_1-t_2 &  0 &  t_1\\[0.2cm]
			0 &  -\frac{1}{64}-3t_1-t_2 &  0 & \frac{5}{32}+4t_1+t_2 & 0 & -\frac{1}{64}-3t_1-t_2 & 0\\[0.2cm]
			t_2 & 0 & \frac{5}{32}+4t_1+t_2 & \frac{1}{2} & \frac{5}{32}+4t_1+t_2 & 0 & t_2\\[0.2cm]
			0 & -\frac{1}{64}-3t_1-t_2 & 0 & \frac{5}{32}+4t_1+t_2 &  0 &  -\frac{1}{64}-3t_1-t_2 & 0\\[0.2cm]
			t_1 &  0 & -\frac{1}{64}-3t_1-t_2 &  0 &  -\frac{1}{64}-3t_1-t_2 & 0 & t_1\\[0.2cm]
			0 &  t_1 &  0 &  t_2 &  0 &  t_1 &  0\end{bmatrix}_{[-3,3]^2},}\ee
	where $t_1,t_2\in\R$ are free parameters. Choose $\{t_1=t_2=0\}$ and $\{t_1=0,\,t_2=\frac{1}{64}\}$, we obtain the following two filters $a,\tta\in l_0(\Z^2)$:
	$${\small a=\frac{1}{64}\begin{bmatrix}0 & -1 & 0 & -1 & 0\\
		-1 & 0 & 10 & 0 & -1\\
		0 & 10 & 32 & 10 & 0\\
		-1 & 0 & 10 & 0 & -1\\
		0 & -1 & 0 & -1 & 0\end{bmatrix}_{[-2,2]^2},\quad\tta=\frac{1}{64}\begin{bmatrix}0 & 0 & 0 & 1 & 0 & 0 & 0\\
		0 & 0 & -2 & 0 & -2 & 0 & 0\\
		0 & -2 & 0 & 11 & 0 & -2 & 0\\
		1 & 0 & 11 & 32 & 11 & 0 & 1\\
		0 & -2 & 0 & 11 & 0 & -2 & 0\\
		0 & 0 & -2 & 0 & -2 & 0 & 0\\
		0 & 0 & 0 & 1 & 0 & 0 & 0\end{bmatrix}_{[-3,3]^2}.}$$
	Choose $n_1=n_2=2$. We now apply Theorem~\ref{thm:df:int} to construct an interpolatory dual $M_{\sqrt{2}}$-framelet from the filters $a$ and $\tta$, with $5$ pairs of interpolatory high-pass filters with symmetry and at least $2$ vanishing moments. With $\dm=M_{\sqrt{2}}$, recall that $\Gamma_{\dm}=\{(0,0)^\tp,\,e_1\}$ is given as in \er{ga:sq2}. The first pair $b_1,\ttb_1\in l_2(\Z^0)$ are defined via \er{b1} and and the second pair $b_2,\ttb_2\in l_0(\Z^2)$ are defined via
	$$\wh{b_2}(\xi)=\frac{1}{\sqrt{2}}-\sqrt{2}e^{-ie_1\cdot\xi}\wh{a^{[e_1,M_{\sqrt{2}}]}}(M_{\sqrt{2}}^\tp\xi),\quad \wh{\ttb_2}(\xi)=\frac{1}{\sqrt{2}}-\sqrt{2}e^{-ie_1\cdot\xi}\wh{\tta^{[e_1,M_{\sqrt{2}}]}}(M_{\sqrt{2}}^\tp\xi),$$
	for all $\xi\in\R^2$. Specifically, we have
	{\small $$b_1=\frac{1}{64}\begin{bmatrix}0 & 1 & 0 & 1 & 0\\
		1 & 0 & -10 & 0 & 1\\
		0 & -10 & 32 & -10 & 0\\
		1 & 0 & -10 & 0 & 1\\
		0 & 1 & 0 & 1 & 0\end{bmatrix}_{[-2,2]^2},\quad\ttb_1=\frac{1}{64}\begin{bmatrix}0 & 0 & 0 & 1 & 0 & 0 & 0\\
		0 & 0 & -2 & 0 & -2 & 0 & 0\\
		0 & -2 & 0 & 11 & 0 & -2 & 0\\
		1 & 0 & 11 & -32 & 11 & 0 & 1\\
		0 & -2 & 0 & 11 & 0 & -2 & 0\\
		0 & 0 & -2 & 0 & -2 & 0 & 0\\
		0 & 0 & 0 & 1 & 0 & 0 & 0\end{bmatrix}_{[-3,3]^2},$$
		$$b_2=\frac{\sqrt{2}}{64}\begin{bmatrix}0 & 1 & 0 & 1 & 0\\
		1 & 0 & -10 & 0 & 1\\
		0 & -10 & 32 & -10 & 0\\
		1 & 0 & -10 & 0 & 1\\
		0 & 1 & 0 & 1 & 0\end{bmatrix}_{[-2,2]^2},\quad\ttb_2=-\frac{\sqrt{2}}{64}\begin{bmatrix}0 & 0 & 0 & 1 & 0 & 0 & 0\\
		0 & 0 & -2 & 0 & -2 & 0 & 0\\
		0 & -2 & 0 & 11 & 0 & -2 & 0\\
		1 & 0 & 11 & -32 & 11 & 0 & 1\\
		0 & -2 & 0 & 11 & 0 & -2 & 0\\
		0 & 0 & -2 & 0 & -2 & 0 & 0\\
		0 & 0 & 0 & 1 & 0 & 0 & 0\end{bmatrix}_{[-3,3]^2}.$$}
	We have $\vmo(b_1)=\vmo(\ttb_1)=4$ and both $b_1,\ttb_1$ are $D_4$-symmetric about $(0,0)^\tp$. Moreover, as $b_2=\sqrt{2}b_1$ and $\ttb_2=-\sqrt{2}\ttb_1$, it follows that $\vmo(b_2)=\vmo(\ttb_2)=4$ and both $b_1,b_2$ are $D_4$-symmetric about $(0,0)^\tp$.
	
	Next, define $h\in l_0(\Z^2)$ via
	$$\wh{h}(\xi)=\frac{1}{2}-2\ol{\wh{a^{[e_1,M_{\sqrt{2}}]}}(\xi)}\wh{\tta^{[e_1,M_{\sqrt{2}}]}}(\xi),\quad\xi\in\R^2.$$
	Since both $a^{[e_1,M_{\sqrt{2}}]}$ and $\tta^{[e_1,M_{\sqrt{2}}]}$ have symmetry type $(D_4,(-1/2,-1/2)^\tp,1)$, the filter $h$ has symmetry type $(D_4,(0,0)^\tp,1)$. Moreover, $\wh{h}(\xi)=\bo(\|\xi\|^4)$ as $\xi\to (0,0)^\tp$. We can construct filters $u_j,\ttu_j\in l_0(\Z^2)$, $j=1,2,3$ such that $\wh{h}=\sum_{j=1}^3\ol{\wh{u_j}}\wh{\ttu_j}$, $\wh{u_j}(\xi)=\bo(\|\xi\|^2)$ and $\wh{\ttu_j}(\xi)=\bo(\|\xi\|^2)$ as $\xi\to (0,0)^\tp$ for all $j=1,2,3$. The rest of the 3 pairs of high-pass filters $b_l,\ttb_l\in l_0(\Z^2)$, $l=3,4,5$ are defined via
	$$\wh{b_l}(\xi):=e^{-ie_1\cdot\xi}\wh{u_{l-2}}(M_{\sqrt{2}}^\tp\xi),\quad \wh{\ttb_l}(\xi):=e^{-ie_1\cdot\xi}\wh{\ttu_{l-2}}(M_{\sqrt{2}}^\tp\xi),\quad\xi\in\R^2,\quad l=3,4,5.$$
	Specifically, the filters $b_l,\ttb_l$, $l=3,4,5$ are given by
	{\small$$b_3=\begin{bmatrix}0 &0 &1\\
		0 &-2 &0\\
		1 &0 &0\end{bmatrix}_{[0,2]\times[-1,1]},\quad \ttb_3=\frac{1}{2048}\begin{bmatrix}0 &0 &1 &0 &0 &0 &0 &0 &0\\
		0 &-2 &0 & -7 & 0 & 0 & 0 & 0 & 0\\
		1 & 0 & 14 & 0 & 0 & 0 & -4 & 0 & 0\\
		0 & -7 & 0 & 0 & 0 & 76 & 0 & 0 & 0\\
		0 & 0 & 0 & 0 & -144 & 0 & 0 & 0 & 0\\
		0 & 0 & 0 & 76 & 0 & 0 & 0 & -7 & 0\\
		0 & 0 & -4 & 0 & 0 & 0 & 14 & 0 & 1\\
		0 & 0 & 0 & 0 & 0 & -7 & 0 & -2 & 0\\
		0 & 0 & 0 & 0 & 0 & 0 & 1 & 0 & 0\end{bmatrix}_{[-3,5]\times[-4,4]},$$	
		$$b_4=\begin{bmatrix}0 &-1 &0\\
		1 & 0 &1\\
		0 & -1 &0\end{bmatrix}_{[0,2]\times[-2,0]},\quad \ttb_4=\frac{1}{2048}\begin{bmatrix}0 & -3 & 0 & 0 & 0 & -3 & 0\\
		3 & 0 & 8 & 0 & 8 & 0 & 3\\
		0 & -8 & 0 & 78 & 0 & -8 & 0\\
		0 & 0 & -78 & 0 & -78 & 0 & 0\\
		0 & -8 & 0 & 78 & 0 & -8 & 0\\
		3 & 0 & 8 & 0 & 8 & 0 & 3\\
		0 & -3 & 0 & 0 & 0 & -3 & 0\end{bmatrix}_{[-2,4]\times[-4,2]},$$
		$$b_5=\begin{bmatrix}1 &0 &0\\
		0 &-2 &0\\
		0 &0 &1\end{bmatrix}_{[-1,1]\times[-2,0]},\quad \ttb_5=\frac{1}{2048}\begin{bmatrix}0 & 0 & 0 & 0 & 0 & 0 & 1 & 0 & 0\\
		0 & 0 & 0 & 0 & 0 & -7 & 0 & -2 & 0\\
		0 & 0 & -4 & 0 & 0 & 0 & 14 & 0 & 1\\
		0 & 0 & 0 & 76 & 0 & 0 & 0 & -7 & 0\\
		0 & 0 & 0 & 0 & -144 & 0 & 0 & 0 & 0\\
		0 & -7 & 0 & 0 & 0 & 76 & 0 & 0 & 0\\
		1 & 0 & 14 & 0 & 0 & 0 & -4 & 0 & 0\\
		0 & -2 & 0 & -7 & 0 & 0 & 0 & 0 & 0\\
		0 & 0 & 1 & 0 & 0 & 0 & 0 & 0 & 0\end{bmatrix}_{[-4,4]\times[-5,3]}.$$}
	We have $\vmo(b_l)=\vmo(\ttb_l)=2$ for all $l=3,4,5$. Moreover, define the following two subgroups of $D_4$:
	$$H_1:=\left\{\pm I_2,\, \pm\begin{bmatrix}0 & 1\\
	1 &0\end{bmatrix}\right\},\quad H_2:=\left\{\pm I_2,\, \pm\begin{bmatrix}1 & 0\\
	0 & -1\end{bmatrix}\right\}.$$
	One can check that $b_3$ and $\ttb_3$ have symmetry type $(H_1,(1,0)^\tp,1)$, $b_4$ and $\ttb_4$ have symmetry type $(H_2,(1,-1)^\tp,1)$, $b_5$ and $\ttb_5$ have symmetry type $(H_1,(0,-1)^\tp,1)$.
	
	Now define $\phi,\tphi$ via \er{ref:tphi} with $\dm=M_{\sqrt{2}}$. By computation, we have $\sm_2(a,M_{\sqrt{2}})\approx 2.4479$ and $\sm_2(\tta,M_{\sqrt{2}})\approx 2.5879$, so by \er{sm:inf:2} we have $\sm_\infty(a,M_{\sqrt{2}})\ge 1.4479$ and $\sm_\infty(\tta,M_{\sqrt{2}})\ge 1.5879$. Therefore, $\phi,\tphi\in L_2(\R^2)$  and are both fundamental. See Figure~\ref{fig:ex1} for the graphs of $\phi$ and $\tphi$. Now define $\psi_l,\tpsi_l$ via \er{ref:tpsi} with $s=5$ and $\dm=M_{\sqrt{2}}$, we then obtain an interpolatory dual $M_{\sqrt{2}}$-framelet $(\{\phi;\psi_1,\dots,\psi_5\},$ $\{\tphi_1;\tpsi_1,\dots,\tpsi_5\})$ in $L_2(\R^2)$ such that $\vmo(\psi_l)=\vmo(\tpsi_l)=4$ for $l=1,2$ and $\vmo(\psi_l)=\vmo(\tpsi_l)=2$ for $l=3,4,5$. Furthermore, $\phi$ and $\tphi$ are $D_4$-symmetric about $(0,0)^\tp$, and using \cite[Proposition 2.1]{han04},the symmetry types of $\psi_l,\tpsi_l$ ($l=1,\dots,5$) are given as follows:
	\begin{itemize}
\item for $l=1,2$, $\psi_l,\tpsi_l$ are $D_4$-symmetric about $(0,0)^\tp$;

	\item $\psi_3$ and $\tpsi_3$ are $H_1$-symmetric about  $c_3:=M_{\sqrt{2}}^{-1}(1,0)^\tp=(1/2,1/2)^\tp$;
		
		\item $\psi_4$ and $\tpsi_4$ are $H_2$-symmetric about  $c_4:=M_{\sqrt{2}}^{-1}(1,-1)^\tp=(0,1)^\tp=e_2$;
		
		\item $\psi_5$ and $\tpsi_5$ are $H_1$-symmetric about  $c_5:=M_{\sqrt{2}}^{-1}(0,-1)^\tp=(-1/2,1/2)^\tp$.
	\end{itemize}
	
	For this particular example, we can further reduce the number of high-pass filters by observing that $b_2=\sqrt{2}b_1$ and $\ttb_2=-\sqrt{2}\ttb_1$. In this case, we have
	$$\ol{\wh{b_1}(\xi)}\wh{\ttb_1}(\xi+2\pi\omega)+\ol{\wh{b_2}(\xi)}\wh{\ttb_2}(\xi+2\pi\omega)=-\ol{\wh{b_1}(\xi)}\wh{\ttb_1}(\xi+2\pi\omega),\quad\forall \xi\in\R^2,\quad \omega\in\Omega_{M_{\sqrt{2}}}.$$
	Therefore, we conclude that $(\{a;b_1,b_3,b_4,b_5\};\{\tta;-\ttb_1,\ttb_3,\ttb_4,\ttb_5\})$ is an interpolatory dual $M_{\sqrt{2}}$-framelet filter bank
	and thus yields the interpolatory dual $M_{\sqrt{2}}$-framelet $(\{\phi;\psi_1,\psi_3,\psi_4,\psi_5\},\{\tphi;-\tpsi_1,\tpsi_3,\tpsi_4,\tpsi_5\})$ in $L_2(\R^2)$.

\begin{figure}[htbp]
\centering
\begin{subfigure}[b]{0.4\textwidth} \includegraphics[width=\textwidth,height=0.8\textwidth]{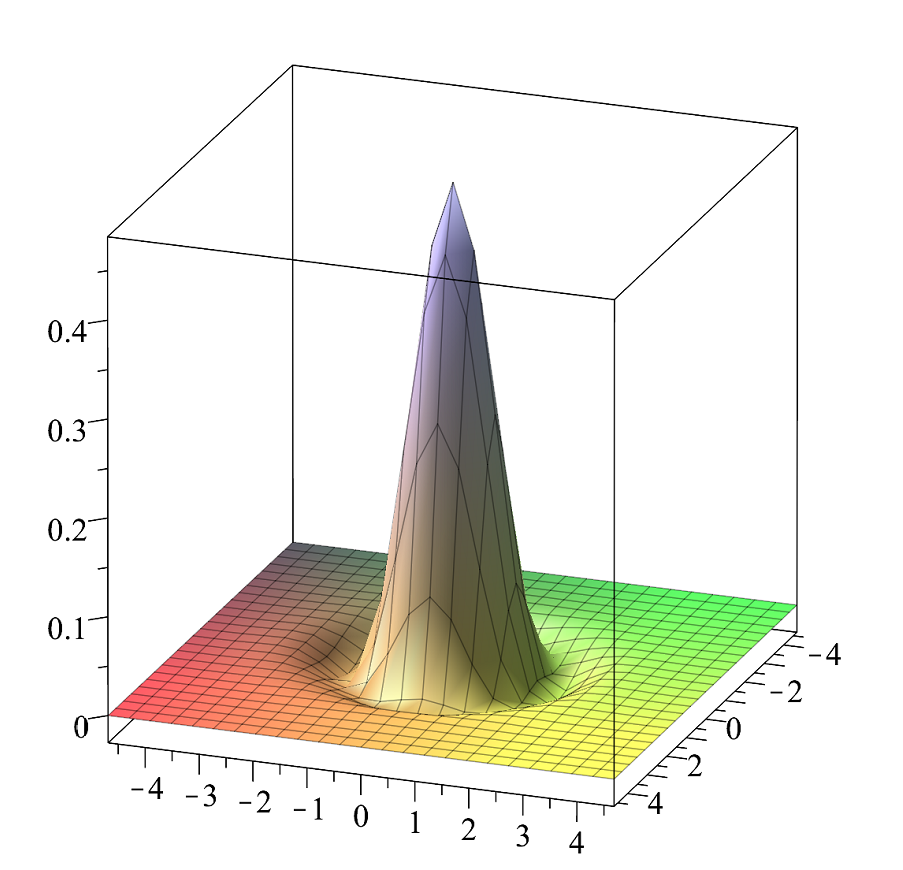}
	\caption{$\phi$}
\end{subfigure}
\begin{subfigure}[b]{0.4\textwidth} \includegraphics[width=\textwidth,height=0.8\textwidth]{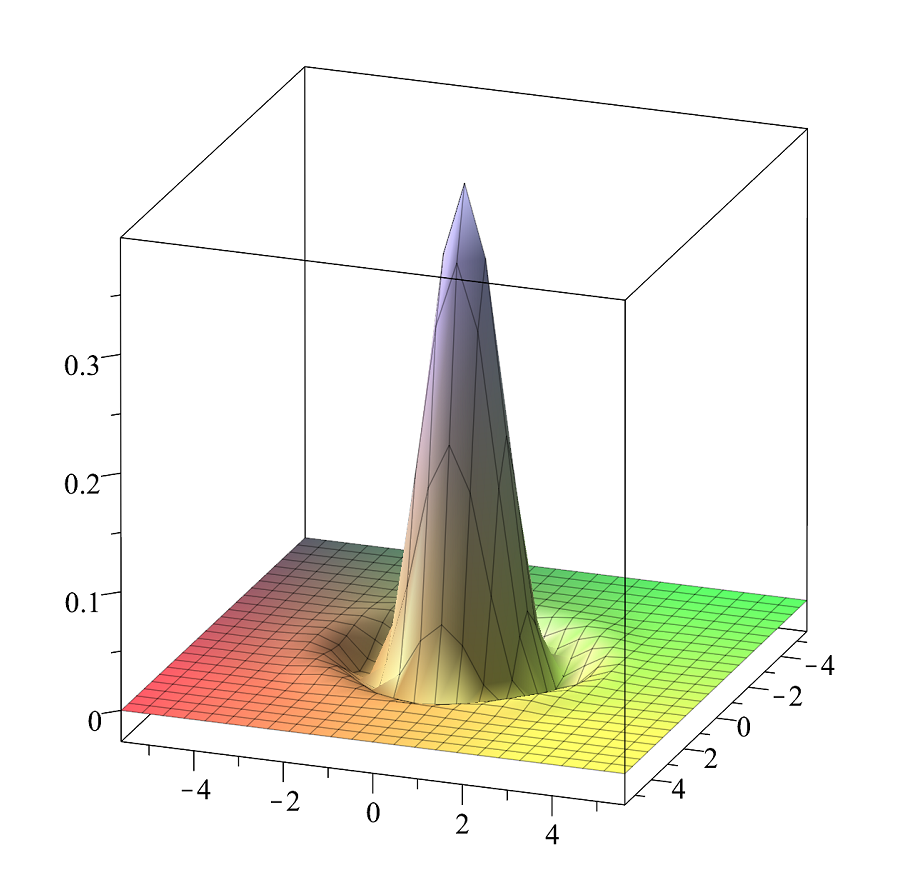}
	\caption{$\tphi$}
\end{subfigure}
\caption{
	The graphs of the interpolatory standard $M_{\sqrt{2}}$-refinable functions $\phi$ and $\tphi$ of the filters $a$ and $\tta$ in Example~\ref{ex1}.
}\label{fig:ex1}
\end{figure}

\end{exmp}

\begin{exmp}\label{ex2}Let $\dm=2I_2$. Parametrize a $2I_2$-interpolatory filter $A\in l_0(\Z^2)$ such that
	\begin{itemize}
		\item $\wh{A}(0,0)=1$ and $\sr(A,2I_2)=4$; 
		
		\item $A$ is $D_6$-symmetric about $(0,0)^\tp$;
	\end{itemize}
	as the following:
	\be\label{A:ex2} {\small A=\begin{bmatrix}0 & 0 & 0 & t & -\frac{1}{64}-t & -\frac{1}{64}-t & t\\[0.2cm]
			0 & 0 & -\frac{1}{64}-t & 0 & \frac{1}{32}+2t & 0 & -\frac{1}{64}-t\\[0.2cm]
			0 & -\frac{1}{64}-t & \frac{1}{32}+2t & \frac{1}{8}-t & \frac{1}{8}-t & \frac{1}{32}+2t & -\frac{1}{64}-t\\[0.2cm]
			t & 0 & \frac{1}{8}-t & \frac{1}{4} & \frac{1}{8}-t & 0 & t\\[0.2cm]
			-\frac{1}{64}-t & \frac{1}{32}+2t & \frac{1}{8}-t & \frac{1}{8}-t & \frac{1}{32}+2t & -\frac{1}{64}-t & 0\\[0.2cm]
			-\frac{1}{64}-t & 0 & \frac{1}{32}+2t & 0 & -\frac{1}{64}-t & 0 & 0\\[0.2cm]
			t & -\frac{1}{64}-t & -\frac{1}{64}-t & t & 0 & 0 & 0\end{bmatrix}_{[-3,3]^2}},\ee
	where $t\in\R$ is a free parameter. Choose $\{t=0\}$ and $\{t=-\frac{1}{64}\}$, we obtain the following two filters $a,\tta\in l_0(\Z^2)$:
	{\small $$a=\frac{1}{64}\begin{bmatrix}0 & 0 & 0 & 0 & -1 & -1 & 0\\
		0 & 0 & -1 & 0 & 2 & 0 & -1\\
		0 & -1 & 2 & 8 & 8 & 2 & -1\\
		0 & 0 & 8 & 16 & 8 & 0 & 0\\
		-1 & 2 & 8 & 8 & 2 & -1 & 0\\
		-1 & 0 & 2 & 0 & -1 & 0 & 0\\
		0 & -1 & -1 & 0 & 0 & 0 & 0\end{bmatrix}_{[-3,3]^2},\quad\tta=\frac{1}{64}\begin{bmatrix}0 & 0 & 0 & -1 & 0 & 0 & -1\\
		0 & 0 & 0 & 0 & 0 & 0 & 0\\
		0 & 0 & 0 & 9 & 9 & 0 & 0\\
		-1 & 0 & 9 & 16 & 9 & 0 & -1\\
		0 & 0 & 9 & 9 & 0 & 0 & 0\\
		0 & 0 & 0 & 0 & 0 & 0 & 0\\
		-1 & 0 & 0 & -1 & 0 & 0 & 0\end{bmatrix}_{[-3,3]^2}.$$}
	Choose $n_1=n_2=2$. We now apply Theorem~\ref{thm:df:int} to construct an interpolatory dual $2I_2$-framelet from the filters $a$ and $\tta$, with $13$ pairs of interpolatory high-pass filters with symmetry and at least $2$ vanishing moments. With $\dm=2I_2$, recall that $\Gamma_{\dm}=\{(0,0)^\tp,\,e_1,\,e_2,\,e_1+e_2\}$ is given as in \er{ga:2id}. Define $b_1,\ttb_1\in l_0(\Z^2)$ via \er{b1} and define $b_l,\ttb_l\in l_0(\Z^2)$ for $l=2,3,4$ via
	$$\wh{b_2}(\xi)=\frac{1}{2}-2e^{-ie_1\cdot\xi}\wh{a^{[e_1,2I_2]}}(2\xi),\quad \wh{\ttb_2}(\xi)=\frac{1}{2}-2e^{-ie_1\cdot\xi}\wh{\tta^{[e_1,2I_2]}}(2\xi),$$
	$$\wh{b_3}(\xi)=\frac{1}{2}-2e^{-ie_2\cdot\xi}\wh{a^{[e_2,2I_2]}}(2\xi),\quad \wh{\ttb_3}(\xi)=\frac{1}{2}-2e^{-ie_2\cdot\xi}\wh{\tta^{[e_2,2I_2]}}(2\xi),$$
	$$\wh{b_4}(\xi)=\frac{1}{2}-2e^{-i(e_1+e_2)\cdot\xi}\wh{a^{[(e_1+e_2),2I_2]}}(2\xi),\quad \wh{\ttb_4}(\xi)=\frac{1}{2}-2e^{-i(e_1+e_2)\cdot\xi}\wh{\tta^{[(e_1+e_2),2I_2]}}(2\xi),$$
	for all $\xi\in\R^2$. Specifically,
	{\small
		$$b_1=\frac{1}{64}\begin{bmatrix}0 & 0 & 0 & 0 & -1 & -1 & 0\\
		0 & 0 & -1 & 0 & 2 & 0 & -1\\
		0 & -1 & 2 & 8 & 8 & 2 & -1\\
		0 & 0 & 8 & -48 & 8 & 0 & 0\\
		-1 & 2 & 8 & 8 & 2 & -1 & 0\\
		-1 & 0 & 2 & 0 & -1 & 0 & 0\\
		0 & -1 & -1 & 0 & 0 & 0 & 0\end{bmatrix}_{[-3,3]^2},\quad\tta=-\frac{1}{64}\begin{bmatrix}0 & 0 & 0 & -1 & 0 & 0 & -1\\
		0 & 0 & 0 & 0 & 0 & 0 & 0\\
		0 & 0 & 0 & 9 & 9 & 0 & 0\\
		-1 & 0 & 9 & -48 & 9 & 0 & -1\\
		0 & 0 & 9 & 9 & 0 & 0 & 0\\
		0 & 0 & 0 & 0 & 0 & 0 & 0\\
		-1 & 0 & 0 & -1 & 0 & 0 & 0\end{bmatrix}_{[-3,3]^2},$$
		$$b_2=\frac{1}{32}\begin{bmatrix}0 & 0 & 1 & 0 & -2 & 0 & 1\\
		0 & 0 & 0 & 0 & 0 & 0 & 0\\
		0 & 0 & -8 & 16 & -8 & 0 & 0\\
		0 & 0 & 0 & 0 & 0 & 0 & 0\\
		1 & 0 & -2 & 0 & 1 & 0 & 0\end{bmatrix}_{[-3,3]\times[-2,2]},\quad \ttb_2=\frac{1}{32}\begin{bmatrix}1 & 0 & -9 & 16 & -9 & 0 & 1
		\end{bmatrix}_{[-3,3]\times\{0\}},$$
		$$b_3=\begin{bmatrix}0 & 0 & 0 & 0 & 1\\
		0 & 0 & 0 & 0 & 0\\
		1 & 0 & -8 & 0 & -2\\
		0 & 0 & 16 & 0 & 0\\
		-2 & 0 & -8 & 0 & 1\\
		0 & 0 & 0 & 0 & 0\\
		1 & 0 & 0 & 0 & 0\end{bmatrix}_{[-2,2]\times[-3,3]},\quad \ttb_3=\begin{bmatrix}1\\
		0\\
		-9\\
		16\\
		-9\\
		0\\
		1\end{bmatrix}_{\{0\}\times[-3,3]},$$
		$$b_4=\begin{bmatrix}0 & 0 & 0 & 0 & 1 & 0 & 0\\
		0 & 0 & 0 & 0 & 0 & 0 &0\\
		0 & 0 & -2 & 0 & -8 & 0 & 1\\
		0 & 0 & 0 & 16 & 0 & 0 & 0\\
		1 & 0 & -8 & 0 & -2 & 0 & 0\\
		0 & 0 & 0 & 0 & 0 & 0 & 0 \\
		0 & 0 & 1 & 0 & 0 & 0 & 0\end{bmatrix}_{[-3,3]^2},\quad \ttb_4=\begin{bmatrix}0 & 0 & 0 & 0 & 0 & 0 & 1\\
		0 & 0 & 0 & 0 & 0 & 0 &0\\
		0 & 0 & 0 & 0 & -9 & 0 & 0\\
		0 & 0 & 0 & 16 & 0 & 0 & 0\\
		0 & 0 & -9 & 0 & 0 & 0 & 0\\
		0 & 0 & 0 & 0 & 0 & 0 & 0 \\
		1 & 0 & 0 & 0 & 0 & 0 & 0\end{bmatrix}_{[-3,3]^2}.$$
	}
	We have $\vmo(b_1)=\vmo(\ttb_1)=4$ and both $b_1,\ttb_1$ have symmetry type $(D_6,(0,0)^\tp,1)$. Using Lemma~\ref{lem:lpm}, we have $\vmo(b_l)=\vmo(\ttb_l)=4$ for $l=2,3,4$. For $j=1,2,3$, define the subgroups $\mathcal{G}_j$ of $D_6$ as in \er{sym:2id}. Using the symmetry types of the coset filters of $a$ and $\tta$, we see that the symmetry types of $b_l$ and $\ttb_l$ are $(\mathcal{G}_{l-1},(0,0)^\tp,1)$ for $l=2,3,4$.
	
	Next, define $h_2, h_3, h_4\in l_0(\Z^2)$ via
	$$\wh{h_2}(\xi)=\frac{1}{4}-4\ol{\wh{a^{[e_1,2I_2]}}(\xi)}\wh{\tta^{[e_1,2I_2]}}(\xi),\quad \wh{h_3}(\xi)=\frac{1}{4}-4\ol{\wh{a^{[e_2,2I_2]}}(\xi)}\wh{\tta^{[e_2,2I_2]}}(\xi),$$
	$$\wh{h_4}(\xi)=\frac{1}{4}-4\ol{\wh{a^{[e_1+e_2,2I_2]}}(\xi)}\wh{\tta^{[e_1+e_2,2I_2]}}(\xi),$$
	for all $\xi\in\R^2$. Using the symmetry of the coset filtes of $a$ and $\tta$, we see that the symmetry type of $h_l$ is $(\mathcal{G}_{l-1},(0,0)^\tp,1)$ for $l=2,3,4$. Moreover, $\wh{h_l}(\xi)=\bo(\|\xi\|^4)$ as $\xi\to (0,0)^\tp$ for $l=2,3,4$. We then find filters $u_{l,j},\ttu_{l,j}\in l_0(\Z^2)$ for $l=2,3,4$ and $j=1,2,3$ such that 
	$$\wh{h_j}(\xi)=\sum_{j=1}^3\ol{\wh{u_{l,j}}(\xi)}\wh{\ttu_{l,j}}(\xi),\quad\forall \xi\in\R^2,\quad l=2,3,4,$$
	and
	$$\wh{u_{l,j}}(\xi)=\bo(\|\xi\|^2),\quad \wh{\ttu_{l,j}}(\xi)=\bo(\|\xi\|^2),\quad\xi\to (0,0)^\tp,\quad\forall l=2,3,4,\quad j=1,2,3.$$
	Specifically, the filters $u_{l,j},\ttu_{l,j}$, $l=2,3,4$, $j=1,2,3$ are given by
	$$\wh{u_{2,1}}(\xi)=(1-e^{i\xi_1})^2,\quad \wh{u_{2,2}}(\xi)=\wh{u_{2,3}}(\xi)=(1-e^{i\xi_1})(1-e^{-i\xi_2}),$$
	$$\wh{\ttu_{2,1}}(\xi)=\frac{1}{1024}(1-e^{i\xi_1})^2[-e^{-i(\xi_1+\xi_2)}-e^{i(\xi_1+\xi_2)}+7(e^{-i\xi_2}+e^{i\xi_2}],$$
	$$ \wh{\ttu_{2,2}}(\xi)=\frac{1}{64}(1-e^{i\xi_1})(1-e^{-i\xi_2}),\quad \wh{\ttu_{2,3}}(\xi)=-\frac{e^{i\xi_1}}{128}(1-e^{-i\xi_1})^2(1+e^{i(\xi_1-\xi_2)}),$$
	$$\wh{u_{3,1}}(\xi)=(1-e^{i\xi_2})^2,\quad \wh{u_{3,2}}(\xi)=(1-e^{i\xi_1})(e^{i\xi_2}-1),\quad \wh{u_{3,3}}(\xi)=(1-e^{i\xi_2})(1-e^{-i\xi_1}),$$
	$$\wh{\ttu_{3,1}}(\xi)=\frac{1}{1024}(1-e^{i\xi_2})^2[-e^{-i(\xi_1+\xi_2)}-e^{i(\xi_1+\xi_2)}+7(e^{-i\xi_1}+e^{i\xi_1}],$$
	$$\wh{\ttu_{3,2}}(\xi)=\frac{1}{64}(1-e^{i\xi_1})(e^{i\xi_2}-1),\quad \wh{u_{3,3}}(\xi)=-\frac{e^{i\xi_2}}{128}(1-e^{-i\xi_2})^2(1+e^{i(\xi_2-\xi_1)}),$$
	$$\wh{u_{4,1}}(\xi)=(1-e^{i\xi_1})(1-e^{-i\xi_1}),\quad \wh{u_{4,2}}(\xi)=(1-e^{i\xi_2})(1-e^{-i\xi_2}),\quad \wh{u_{4,3}}(\xi)=(1-e^{i\xi_1})(1-e^{i\xi_2}),$$
	$$\wh{\ttu_{4,1}}(\xi)=\frac{e^{i(\xi_1+\xi_2)}}{1024}(1-e^{-i(\xi_1+\xi_2)})^2[e^{-i\xi_1}+e^{-\xi_1}-14],$$
	$$\wh{u_{4,2}}(\xi)=\frac{e^{i(\xi_1+\xi_2)}}{1024}(1-e^{-i(\xi_1+\xi_2)})^2[7(e^{-i\xi_2}+e^{-\xi_2})-6(e^{-i\xi_1}+e^{i\xi_1}+e^{i(\xi_1-\xi_2)}+e^{i(\xi_2-\xi_1)})-2],$$
	$$\wh{u_{4,3}}(\xi)=\frac{e^{i(\xi_1+\xi_2)}}{1024}(1-e^{-i(\xi_1+\xi_2)})^2[6(e^{i(\xi_1-\xi_2)}+e^{2i\xi_2})-(e^{-i\xi_1}+e^{-i\xi_2}+e^{i(2\xi_1+2\xi_2)}+e^{i(\xi_1+2\xi_1)})],$$
	for all $\xi=(\xi_1,\xi_2)^\tp\in\R^2$. Note that $u_{l,j}$ has symmetry for all $l=2,3,4$ and $j=1,2,3$. Define the rest of the high-pass filters $b_l,\ttb_l\in l_0(\Z^2)$, $l=5,\dots,13$ via
	$$\wh{b_{4+j}}(\xi):=e^{-ie_1\cdot\xi}\wh{u_{2,j}}(2\xi),\quad \wh{\ttb_{4+j}}(\xi):=e^{-ie_1\cdot\xi}\wh{\ttu_{2,j}}(2\xi),$$
	$$\wh{b_{7+j}}(\xi):=e^{-ie_2\cdot\xi}\wh{u_{3,j}}(2\xi),\quad \wh{\ttb_{7+j}}(\xi):=e^{-ie_2\cdot\xi}\wh{\ttu_{3,j}}(2\xi),$$
	$$\wh{b_{10+j}}(\xi):=e^{-i(e_1+e_2)\cdot\xi}\wh{u_{4,j}}(2\xi),\quad \wh{\ttb_{10+j}}(\xi):=e^{-i(e_1+e_2)\cdot\xi}\wh{\ttu_{4,j}}(2\xi),$$
	for all $\xi\in\R^2$ and for $j=1,2,3$. For simplicity, we do not present these filters here. By our construction, we have $\vmo(b_l)=\vmo(\ttb_l)=2$ for all $l=5,\dots,13$, and the symmetry types of these high-pass filters are as follows:
	\begin{center}\begin{tabular}{|c|c|c|c|}
			\hline 
			Filter $b_l$ & Symmetry type of $b_l$ & Filter $\ttb_l$ & Symmetry type of $\ttb_l$ \\ 
			\hline 
			$b_5$ & $(\cG_1,(-1,0)^\tp,1)$ & $\ttb_5$ & $(\{\pm I_2\},(-1,0)^\tp,1)$ \\ 
			\hline 
			$b_6$ & $(\cG_3,(0,1)^\tp,1)$ & $\ttb_6$ & $(\cG_3,(0,1)^\tp,1)$ \\ 
			\hline 
			$b_7$ & $(\cG_3,(0,1)^\tp,1)$ & $\ttb_7$ & $(\{\pm I_2\},(0,1)^\tp,1)$ \\ 
			\hline 
			$b_8$ & $(\cG_1,(0,0)^\tp,1)$ & $\ttb_8$ & $(\{\pm I_2\},(0,0)^\tp,1)$ \\ 
			\hline 
			$b_9$ & $(\cG_3,(-1,0)^\tp,1)$ & $\ttb_9$ & $(\cG_3,(-1,0)^\tp,1)$ \\ 
			\hline 
			$b_{10}$ & $(\cG_3,(1,0)^\tp,1)$ & $\ttb_{10}$ & $(\{\pm I_2\},(1,0)^\tp,1)$ \\ 
			\hline 
			$b_{11}$ & $(\cG_1,(1,1)^\tp,1)$ & $\ttb_{11}$ & $(\{\pm I_2\},(1,1)^\tp,1)$ \\ 
			\hline 
			$b_{12}$ & $(\cG_2,(1,1)^\tp,1)$ & $\ttb_{12}$ & $(\{\pm I_2\},(1,1)^\tp,1)$ \\ 
			\hline 
			$b_{13}$ & $(\cG_3,(0,0)^\tp,1)$ & $\ttb_{13}$ & $(\{\pm I_2\},(0,0)^\tp,1)$ \\ 
			\hline 
	\end{tabular} \end{center}
	Now define $\phi,\tphi$ via \er{ref:tphi} with $\dm=2I_2$. By computation, we have $\sm_2(a,2I_2)\approx 2.4408$ and $\sm_2(\tta,2I_2)\approx 1.7658$, so $\sm_\infty(a,2I_2)\ge 1.4408$ and $\sm_\infty(\tta,2I_2)\ge 0.7658$. Therefore, $\phi,\tphi\in L_2(\R^2)$  and are both fundamental. See Figure~\ref{fig:ex2} for the graphs of $\phi$ and $\tphi$. Now define $\psi_l,\tpsi_l$ via \er{ref:tpsi} with $s=13$ and $\dm=2I_2$, we then obtain an interpolatory dual $2I_2$-framelet $(\{\phi;\psi_1,\dots,\psi_{13}\},$ $\{\tphi_1;\tpsi_1,\dots,\tpsi_{13}\})$ in $L_2(\R^2)$ such that $\vmo(\psi_l)=\vmo(\tpsi_l)=4$ for $l=1,2,3,4$ and $\vmo(\psi_l)=\vmo(\tpsi_l)=2$ for $l=5,\dots,13$. Furthermore, $\phi$ and $\tphi$ are $D_6$-symmetric about $(0,0)^\tp$, the symmetry types of $\psi_l,\tpsi_l$ $l=1,\dots,13$ can be obtained by using \cite[Proposition 2.1]{han04}, and we omit the details here for the simplicity of presentation.
\begin{figure}[htbp]
\centering
\begin{subfigure}[b]{0.4\textwidth} \includegraphics[width=\textwidth,height=0.8\textwidth]{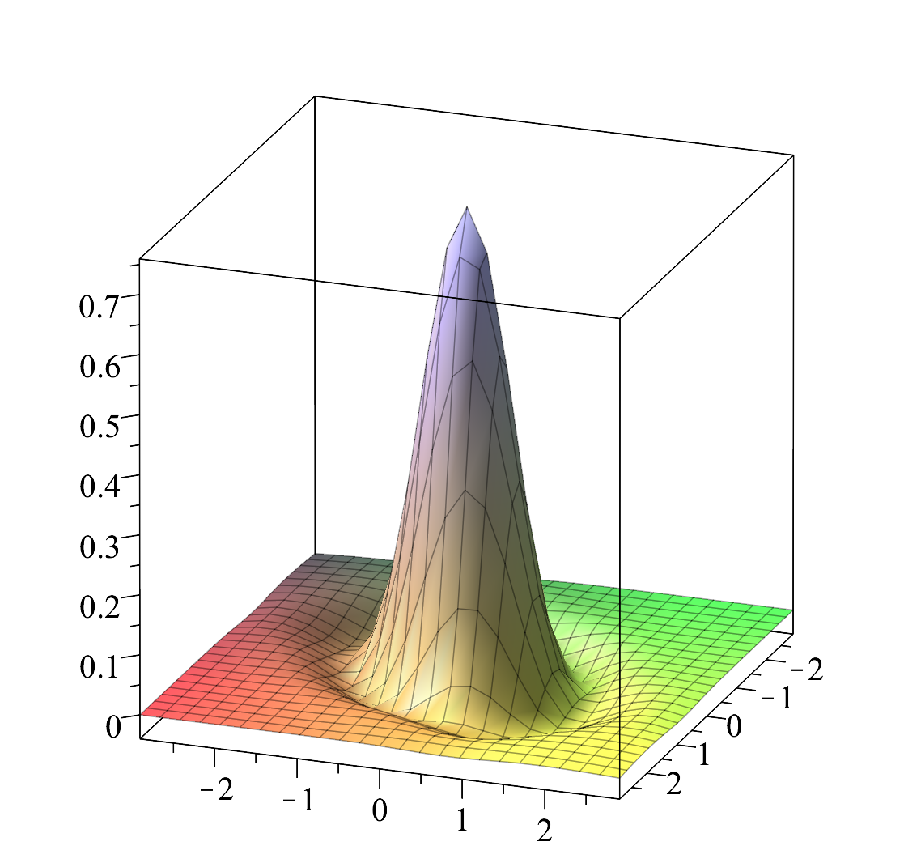}
	\caption{$\phi$}
\end{subfigure}
\begin{subfigure}[b]{0.4\textwidth} \includegraphics[width=\textwidth,height=0.8\textwidth]{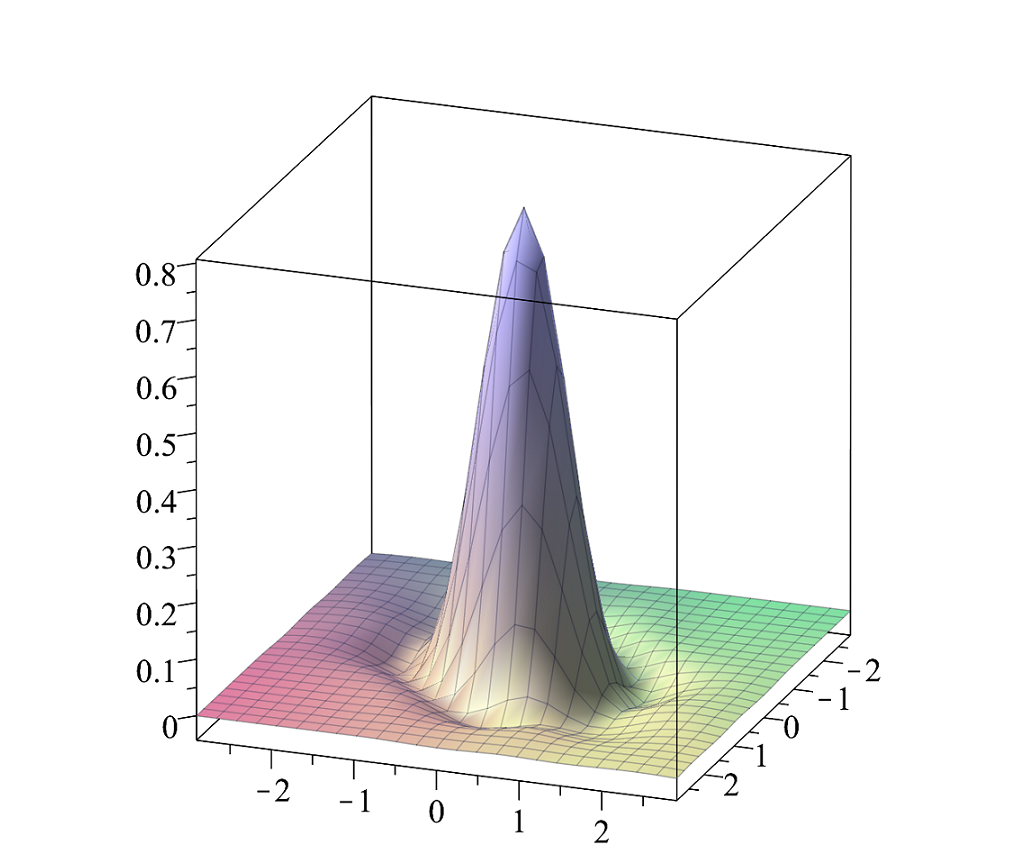}
	\caption{$\tphi$}
\end{subfigure}
\caption{
	The graphs of the interpolatory standard $2I_2$-refinable functions $\phi$ and $\tphi$ of the filters $a$ and $\tta$ in Example~\ref{ex2}.
}\label{fig:ex2}
\end{figure}
	
\end{exmp}

\begin{exmp}\label{ex3}Let $\dm=M_{\sqrt{3}}=\begin{bmatrix}1 & -2\\
	2 &-1\end{bmatrix}$. Let $a\in l_0(\Z^2)$ be given by
	{\small $$a=\frac{1}{243}\begin{bmatrix}0 & 0 & 0 & 0 & -2 & -2 & 0\\
		0 & 0 & -2 & -1 & 0 & -1 & -2\\
		0 & -2 & 0 & 32 & 32 & 0 & -2\\
		0 & -1 & 32 & 81 & 32 & -1 & 0\\
		-2 & 0 & 32 & 32 & 0 & -2 & 0\\
		-2 & -1 & 0 & -1 & -2 & 0 & 0\\
		0 & -2 & -2 & 0 & 0 & 0 & 0\end{bmatrix}_{[-3,3]^2}.$$}
	The filter $a$ satisfies $\wh{a}(0,0)=1$, $\sr(A,M_{\sqrt{3}})=4$, and is $D_6$-symmetric about $(0,0)^\tp$. Choose $m=2$. We apply Corollary~\ref{cor:qtf:int} to construct an interpolatory quasi-tight $M_{\sqrt{3}}$-framelet from $a$. With $\dm=M_{\sqrt{3}}$, recall that $\Gamma_{\dm}=\{(0,0)^\tp,\,-e_1,\,e_2\}$ is given as in \er{ga:sq3}. Define $b_1,b_2,b_3\in\l_0(\Z^2)$ and $\eps_1,\eps_2,\eps_3\in\{-1,1\}$via
	$$\wh{b_1}(\xi)=\wh{a}(\xi)-1,\quad \eps_1=-1,$$
	$$\wh{b_2}(\xi)=\frac{1}{\sqrt{3}}-\sqrt{3}\wh{a^{[-e_1,\dm]}}(M_{\sqrt{3}}^\tp\xi),\quad \eps_2=1,$$
	$$\wh{b_3}(\xi)=\frac{1}{\sqrt{3}}-\sqrt{3}\wh{a^{[e_2,\dm]}}(M_{\sqrt{3}}^\tp\xi),\quad \eps_2=1,$$
	for all $\xi\in\R^2$. Specifically, we have
	{\small $$b_1=\frac{1}{243}\begin{bmatrix}0 & 0 & 0 & 0 & 2 & 2 & 0\\
		0 & 0 & 2 & 1 & 0 & 1 & 2\\
		0 & 2 & 0 & -32 & -32 & 0 & 2\\
		0 & 1 & -32 & 162 & -32 & 1 & 0\\
		2 & 0 & -32 & -32 & 0 & 2 & 0\\
		2 & 1 & 0 & 1 & 2 & 0 & 0\\
		0 & 2 & 2 & 0 & 0 & 0 & 0\end{bmatrix}_{[-3,3]^2},\quad b_2=\frac{\sqrt{3}}{243}\begin{bmatrix}0 & 0 & 0 & 0 & 0 & 2 & 0\\
		0 & 0 & 0 & 1 & 0 & 0 & 2\\
		0 & 2 & 0 & 0 & -32 & 0 & 0\\
		0 & 0 & -32 & 81 & 0 & 1 & 0\\
		2 & 0 & 0 & -32 & 0 & 0 & 0\\
		0 & 1 & 0 & 0 & 2 & 0 & 0\\
		0 & 0 & 2 & 0 & 0 & 0 & 0\end{bmatrix}_{[-3,3]^2},$$
		$$b_3=\frac{\sqrt{3}}{243}\begin{bmatrix}0 & 0 & 0 & 0 & 2 & 0 & 0\\
		0 & 0 & 2 & 0 & 0 & 1 & 0\\
		0 & 0 & 0 & -32 & 0 & 0 & 2\\
		0 & 1 & 0 & 81 & -32 & 0 & 0\\
		0 & 0 & -32 & 0 & 0 & 2 & 0\\
		2 & 0 & 0 & 1 & 0 & 0 & 0\\
		0 & 2 & 0 & 0 & 0 & 0 & 0\end{bmatrix}_{[-3,3]^2}.$$}
	We have $\vmo(b_1)=\vmo(b_2)=\vmo(b_3)=4$. The filter $b_1$ has symmetry type $(D_6,(0,0)^\tp,1)$ and the filters $b_2,b_3$ have symmetry type $(\mathcal{\tilde{H}},(0,0)^\tp,1)$, where 
$$\mathcal{\tilde{H}}:=\dm\mathcal{H}\dm^{-1}=\left\{I_2,\,
\begin{bmatrix}0 & -1\\
		1 &-1\end{bmatrix},\,
\begin{bmatrix}-1 & 1\\
		-1 &0\end{bmatrix},\,
\begin{bmatrix}0 & 1\\
		1 &0\end{bmatrix},\,
\begin{bmatrix}0 & -1\\
		1 &-1\end{bmatrix},\,
\begin{bmatrix}-1 & 0\\
		-1 & 1\end{bmatrix}\right\},$$ 
and $\mathcal{H}$ is the subgroup of $D_6$ as in \er{sym:sq3:1}.
	
	Next, define $h_2,h_3\in l_0(\Z^2)$ via
	$$\wh{h_2}(\xi)=\frac{1}{3}-3|\wh{a^{[-e_1,\dm]}}(\xi)|^2,\quad \wh{h_3}(\xi)=\frac{1}{3}-3|\wh{a^{[e_2,\dm]}}(\xi)|^2,\quad\forall\xi\in\R^2.$$
	Using the symmetry of $a$, we see that $h_2,h_3$ is $\mathcal{H}$-symmetric about $(0,0)^\tp$. In fact, for this particular example, $h_2,h_3$ have symmetry type $(D_6,(0,0)^\tp,1)$. Moreover, $\wh{h_l}(\xi)=\bo(\|\xi\|^4)$ as $\xi\to (0,0)^\tp$ for $l=2,3$. We then find filters $u_{2,j},u_{3,j}\in l_0(\Z^2)$ and $\eps_{2,j},\eps_{3,j}\in\{-1,1\}$ for $j=1,\dots,9$ sicj tjat
	$$\wh{h_l}(\xi)=\sum_{j=1}^9\eps_{l,j}|\wh{u_{l,j}}(\xi)|^2,\quad\forall \xi\in\R^2,\quad l=2,3.$$
	Specifically, $u_{2,j},u_{3,j}$ and $\eps_{2,j},\eps_{3,j}$ are given by
	{\footnotesize $$u_{2,1}=\frac{1}{78732}\begin{bmatrix}-80 & 21083 & -42006 & 21083 & -80	\end{bmatrix}_{[-1,3]\times\{0\}},\quad u_{2,2}=\frac{1}{78732}\begin{bmatrix}-80\\
		21083\\
		-42006\\
		21083\\
		-80\end{bmatrix}_{\{0\}\times[-1,3]},$$
		$$u_{2,3}=\frac{1}{39366}\begin{bmatrix}0 & -12 & 24 & -16 & 8 & -4\\
		0 & 267 & -534 & 257 & 20 & -10\\
		0 & -543 & -37692 & 37631 & 620 & -16\\
		-16 & 620 & 37631 & -37692 & -543 & 0\\
		-10 & 20 & 257 & -534 & 267 & 0\\
		-4 & 8 & -16 & 24 & -12 & 0\end{bmatrix}_{[-3,2]^2},$$
		$$u_{2,4}=\frac{1}{13122}\begin{bmatrix}4 & -8 & 0 & -6\\
		-89 & 180 & 0 & 29\\
		183 & -13458 & 13286 & -121\\
		-121 & 13286 & -13458 & 183\\
		29 & 0 & 180 & -89\\
		-6 & 0 & -8 & 4\end{bmatrix}_{[-2,1]\times[-3,2]}, \quad u_{2,5}=\frac{2\sqrt{357}}{27}\begin{bmatrix}-1 & 1\\1 &-1
		\end{bmatrix}_{[0,1]^2},$$
		$$u_{2,6}=\frac{1}{78732}\begin{bmatrix}80 & 18283 & -36726 & 18283 & 80\end{bmatrix}_{[-1,3]\times\{0\}},\quad u_{2,7}=\frac{1}{78732}\begin{bmatrix}80\\
		18283\\
		-36726\\
		18283\\
		80\end{bmatrix}_{\{0\}\times[-1,3]},$$
		
		$$u_{2,8}=\frac{1}{39366}\begin{bmatrix}0 & -12 & 24 & -16 & 8 & -4\\
		0 & 267 & -534 & 257 & 20 & -10\\
		0 & -543 & 1674 & -1735 & 620 & -16\\
		-16 & 620 & -1735 & 1674 & -543 & 0\\
		-10 & 20 & 257 & -534 & 267 & 0\\
		-4 & 8 & -16 & 24 & -12 & 0\end{bmatrix}_{[-3,2]^2},\quad u_{2,9}=\frac{1}{13122}\begin{bmatrix}4 & -8 & 0 & -6\\
		-89 & 180 & 0 & 29\\
		183 & -336 & 164 & -121\\
		-121 & 164 & -336 & 183\\
		29 & 0 & 180 & -89\\
		-6 & 0 & -8 & 4\end{bmatrix}_{[-2,1]\times[-3,2]},$$
		
		$$u_{3,1}=u_{2,1},\quad u_{3,2}=u_{2,2},\quad u_{3,5}=u_{2,5},\quad u_{3,6}=u_{2,6},\quad u_{3,7}=u_{2,7},$$
		
		$$u_{3,3}=\frac{1}{39366}\begin{bmatrix}0 & -48 & 78 & -16 & -10 & -4\\
		0 & 441 & -795 & 257 & 107 & -10\\
		0 & -732 & -37332 & 37478 & 602 & -16\\
		-16 & 602 & 37478 & -37332 & -732 & 0\\
		-10 & 107 & 257 & -795 & 441 & 0\\
		-4 & -10 & -16 & 78 & -48 & 0\end{bmatrix}_{[-3,2]^2},$$ $$u_{3,4}=\frac{1}{13122}\begin{bmatrix}16 & -26 & 0 & 0\\
		-147 & 267 & 0 & 0\\
		246 & -13578 & 13337 & -115\\
		-115 & 13337 & -13578 & 246\\
		0 & 0 & 267 & -147\\
		0 & 0 & -26 & 16\end{bmatrix}_{[-2,1]\times[-3,2]}$$
		
		$$u_{3,8}=\frac{1}{39366}\begin{bmatrix}0 & -48 & 78 & -16 & -10 & -4\\
		0 & 441 & -795 & 257 & 107 & -10\\
		0 & -732 & 2034 & -1888 & 602 & -16\\
		-16 & 602 & -1888 & 2034 & -732 & 0\\
		-10 & 107 & 257 & -795 & 441 & 0\\
		-4 & -10 & -16 & 78 & -48 & 0\end{bmatrix}_{[-3,3]^2},\quad u_{3,9}=\frac{1}{13122}\begin{bmatrix}16 & -26 & 0 & 0\\
		-147 & 267 & 0 & 0\\
		246 & -456 & 215 & -115\\
		-115 & 215 & -456 & 246\\
		0 & 0 & 267 & -147\\
		0 & 0 & -26 & 16\end{bmatrix}_{[-2,1]\times[-3,2]},$$
	}
	$$\eps_{l,j}=1,\quad l=2,3,\quad j=1,2,3,4,$$
	$$\eps_{l,j}=-1,\quad l=2,3,\quad j=5,6,7,8,9.$$
	Note that $u_{l,j}$ have symmetry for all $l=2,3$ and $j=1,\dots,9$. Define $\eps_{l}\in\{-1,1\}$ and $b_l\in l_0(\Z^2)$ for $l=4,\dots,21$ via
	$$\eps_3+j:=\eps_{2,j},\quad \eps_{12+j}:=\eps_{3,j},\quad \wh{b_{3+j}}(\xi):=e^{ie_1\cdot\xi}\wh{u_{2,j}}(M_{\sqrt{3}}^\tp\xi),\quad \wh{b_{12+j}}(\xi):=e^{-ie_2\cdot\xi}\wh{u_{3,j}}(M_{\sqrt{3}}^\tp\xi),$$
for all $\xi\in\R^2,$ $j=1,\dots,9$	Then $\{a;b_1,\dots,b_{21}\}_{\eps_1,\dots,\eps_{21}}$ forms an interpolatory quasi-tight $M_{\sqrt{3}}$-framelet filter bank. Furthermore, $\vmo(b_l)=2$ for all $l=4,\dots,21$, and the symmetry types of the filters $b_4,\dots,b_{21}$ are given as follows:
	\begin{center}\begin{tabular}{|c|c|c|c|}
			\hline 
			Filter $b_l$ & Symmetry type of $b_l$ & Filter $b_l$ & Symmetry type of $b_l$ \\ 
			\hline 
			$b_4$ & $(\cG_1,(0,2)^\tp,1)$ & $b_{13}$ & $(\cG_1,(1,3)^\tp,1)$  \\ 
			\hline 
			$b_5$ & $(\cG_2,(-3,-1)^\tp,1)$ & $b_{14}$ & $(\cG_2,(-2,0)^\tp,1)$ \\ 
			\hline 
			$b_6$ & $(\{\pm I_2,(-1/2,-1/2)^\tp,1)$ & $b_{15}$ & $(\{\pm I_2,(1/2,1/2)^\tp,1)$ \\ 
			\hline 
			$b_7$ & $(\{\pm I_2,(-1/2,-1/2)^\tp,1)$ & $b_{16}$ & $(\{\pm I_2,(1/2,1/2)^\tp,1)$ \\ 
			\hline 
			$b_8$ & $(\cG_3,(-3/2,1/2)^\tp,1)$ & $b_{17}$ & $(\cG_3,(-1/2,3/2)^\tp,1)$ \\ 
			\hline 
			$b_9$ & $(\cG_1,(0,2)^\tp,1)$ & $b_{18}$ & $(\cG_1,(1,3)^\tp,1)$ \\ 
			\hline 
			$b_{10}$ & $(\cG_2,(-3,-1)^\tp,1)$ & $b_{19}$ & $(\cG_2,(-2,0)^\tp,1)$ \\ 
			\hline 
			$b_{11}$ & $(\{\pm I_2,(-1/2,-1/2)^\tp,1)$ & $b_{20}$ &$(\{\pm I_2,(1/2,1/2)^\tp,1)$  \\ 
			\hline 
			$b_{12}$ & $(\{\pm I_2,(-1/2,-1/2)^\tp,1)$ & $b_{21}$ & $(\{\pm I_2,(1/2,1/2)^\tp,1)$ \\ 
			\hline 
	\end{tabular} \end{center}
	where $\cG_{1},\cG_2,\cG_3$ are the subgroups of $D_6$ that are given in \er{sym:2id}. Now define $\phi$ via \er{ref:phi} with $\dm=M_{\sqrt{3}}$. By computation, we have $\sm_2(a,M_{\sqrt{3}})\approx 2.52996$,  so $\sm_\infty(a,M_{\sqrt{3}})\ge 1.52996$. Therefore, $\phi\in L_2(\R^2)$ and is fundamental. See Figure~\ref{fig:ex3} for the graph of  $\phi$. Now define $\psi_l$ via \er{ref:tpsi} with $s=21$ and $\dm=M_{\sqrt{3}}$, we then obtain an interpolatory quasi-tight $M_{\sqrt{3}}$-framelet $\{\phi;\psi_1,\dots,\psi_{21}\}_{\eps_1,\dots,\eps_{21}}$ in $L_2(\R^2)$ such that $\vmo(\psi_l)=4$ for $l=1,2,3$ and $\vmo(\psi_l)=2$ for $l=4,\dots,21$. Furthermore, $\phi$ is $D_6$-symmetric about $(0,0)^\tp$, the symmetry types of $\psi_1,\dots\psi_{21}$ can be obtained by using \cite[Proposition 2.1]{han04} and we omit the details here for the simplicity of presentation.

\begin{figure}\centering
	\includegraphics[width=0.4\linewidth] {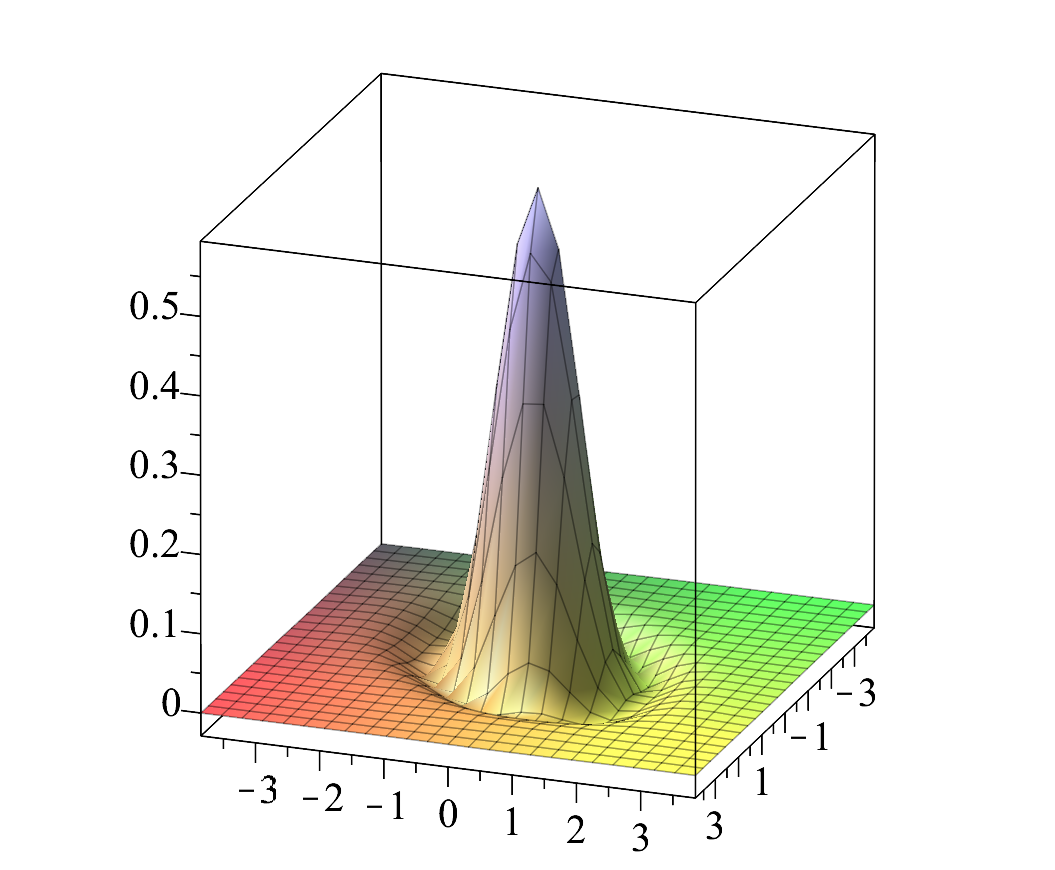}
\caption{The graph of the interpolatory standard $M_{\sqrt{3}}$-refinable function $\phi$ of the filter $a$ in Example~\ref{ex3}.}\label{fig:ex3}
\end{figure}
	
\end{exmp}

\end{document}